\documentclass[oneside,english]{amsart}
\usepackage[T1]{fontenc}
\usepackage[latin9]{inputenc}
\usepackage{amsthm,amsmath,amssymb,color,graphicx,babel}
\usepackage{float}

\newtheorem{theorem}{Theorem}
\newtheorem{lemma}[theorem]{Lemma}
\newtheorem{corollary}[theorem]{Corollary}

\newtheorem{proposition}[theorem]{Proposition}
\theoremstyle{definition}
\newtheorem{definition}[theorem]{Definition}
\newtheorem{remark}[theorem]{Remark}

\newcommand{\mW}{{\mathcal W}}

\newcommand{\Aut}{\hbox{\rm Aut}}

\newcommand{\inv}{^{-1}}
\renewcommand{\mod}{\hbox{mod}\, }

\newcommand{\ZZ}{\mathbb{Z}}

\newcommand{\Cov}{\hbox{\rm Cov}}

\newcommand{\beg}{\mathop{{\rm beg}}}
\renewcommand{\inv}{\mathop{{\rm inv}}}

\newcommand{\D}{{\rm D}}
\newcommand{\V}{{\rm V}}

\renewcommand{\mod}{\hbox{\rm{mod }}}

\makeatletter
\numberwithin{equation}{section}
\numberwithin{figure}{section}

\makeatother

\usepackage{babel}

\title[]{Classification of cubic vertex-transitive tricirculants}
\author{Primo\v{z} Poto\v{c}nik}
\author{Micael Toledo}

\address{Primo\v{z} Poto\v{c}nik, Faculty of Mathematics and Physics, University of Ljubljana, Jadranska 21, SI-1000 Ljubljana, Slovenia.\newline
\indent Also affiliated with: Institute of Mathematics, Physics and Mechanics, Jadranska 19, SI-1000 Ljubljana, Slovenia.
}
\email{primoz.potocnik@fmf.uni-lj.si}

\address{Micael Toledo, Institute of Mathematics, Physics and Mechanics, Jadranska 19, SI-1000.\newline
 Also affiliated with: University of Primorska, Faculty of Mathematics, Natural Sciences and Information Technologies, Glagolja\v{s}ka 8, SI-6000 Koper, Slovenia.}
\email{micael.toledo@imfm.uni-lj.si}	

\begin{document}

\begin{abstract}
A finite graph is called a tricirculant if admits a cyclic group of automorphism which has
precisely three orbits on the vertex-set of the graph, all of equal size.
We classify all finite connected cubic vertex-transitive tricirculants.
We show that except for some small exceptions of order less than $54$,
each of these graphs is either a prism of order $6k$ with $k$ odd, a M\"obius ladder,
or it falls into one of two infinite families,
 each family containing one graph for every order of the form $6k$ with $k$ odd.
\end{abstract}
\maketitle

\section{Introduction}

All graphs in this paper are finite. A connected graph $\Gamma$ admitting a cyclic group of automorphisms $H$ having $k$ orbits of vertices of equal size larger than $1$ is called a {\em $k$-circulant} and a generator of $H$ is then called a {\em $k$-circulant automorphism} of $\Gamma$. In particular, $1$-, $2$- and $3$-circulants are generally called circulants, bicirculants and tricirculants, respectively. A graph is called cubic if it is connected and regular of valence $3$.

A famous and longstanding {\em polycirculant conjecture} claims that every vertex-transitive graph and digraph is a $k$-circulant for some $k$ (see \cite{poly2,poly1}). It is particularly interesting to study conditions under which vertex-transitive graphs admit $k$-circulant automorphisms of large order (and thus relatively small $k$); see, for example, \cite{cubicpoly1,cubicpoly2}.
 On the other hand, existence of a $k$-circulant automorphism with small $k$ often restricts the structure of a vertex-transitive graph to the extent that allows a complete classification. There is a  series of classification  results about
cubic arc-transitive $k$-circulant for $k \le 5$ (see \cite{FK,symtric}). 
In particular, cubic arc-transitive tricirculants where completely classified in \cite{symtric}, where it was shown that only $4$ such graphs exist: $K_{3,3}$, the Pappus Graph, Tutte's $8$-cage and a graph on $54$ vertices. 
This work culminated in a beautiful paper \cite{kmulti}, where it was shown that for all square-free values of $k$ coprime to $6$ there exist only finitely many cubic arc-transitive $k$-circulants.

In this paper we widen the focus to a much wider and structurally richer class of cubic 
vertex-transitive graphs that are not necessarily  
arc-transitive. It is know that every cubic vertex-transitive graph has a $k$-circulant automorphisms \cite{MarSca}
and that no fixed $k$ exists such that every cubic vertex-transitive graph is a $k$-circulant \cite{cubicpoly2}.
As for the classification results, it is clear that a cubiv graph is a vertex-transitive $1$-circulant if and only if it is
a cubic Cayley graph on a cyclic group. Furthermore, vertex-transitive cubic  bicirculants were classified in \cite{bic}.
The goal of this paper is to provide a complete classification of cubic, vertex-transitive tricirculants.

The following theorem and remarks are a brief summary of the contents of this paper.

\begin{theorem}
Let $\Gamma$ be a cubic vertex-transitive tricirculant. Then the order of $\Gamma$ is $6k$ 
for some positive integer $k$ and one of the following holds:
\begin{enumerate}
\item $k$ is odd and $\Gamma$ is isomorphic to one of the graphs $\textrm{X}(k)$ or $\textrm{Y}(k)$, described in Section~\ref{type1} (Definition~\ref{def:X}) and Section~\ref{Type2} (Definition~\ref{def:Y}), respectively.
\item $\Gamma$ is a prism and $k$ is odd, or $\Gamma$ is a M\"obius ladder. 
\item the order of $\Gamma$ is less than $54$ and $\Gamma$ is one of the twenty graph described in Section~\ref{sec:small}.
\end{enumerate}
\end{theorem}

\begin{remark}
Prisms and M\"obius Ladders are circulant graphs. Meanwhile $\textrm{X}(k)$ and $\textrm{Y}(k)$ are bicirculants if $k$ is odd and $3 \nmid k$. In this case, $\textrm{X}(k)$
is isomorphic to the Generalized Petersen graph $\textrm{GP}(3k,k+(-1)^\alpha)$ where $\alpha \in \{1,2\}$ and $\alpha \equiv k$ (\mod$3$).
\end{remark}

\begin{remark}
$\textrm{Y}(k)$ can be embedded on the torus yielding the toroidal map $\{6,3\}_{\alpha,3}$ where $\alpha = \frac{1}{2}(k-1)$ (see \cite{torus} for the definitions pertaining to the maps on the torus).
\end{remark}

\section{Diagrams}

Even though this paper is about simple graphs, it will be 
 convenient in the course of the analysis 
 to have a slightly more general definition, which allows the graphs 
to have loops, parallel edges and semi-edges.
To distinguish between a simple graph (which we will refer to simply as a graph) and these more
general objects, that we shall call {\it pregraphs}.
In what follows, we briefly introduce the concept of a pregraph  and
refer the reader to \cite{covers,elabcovers} for more detailed explanation.

A {\em pregraph} is an ordered $4$-tuple $(V,D; \beg,\inv)$ where
$D$ and $V \neq \emptyset$ are disjoint finite sets of {\em darts}
and {\em vertices}, respectively, $\beg: D \to V$ is a mapping
which assigns to each dart $x$ its {\em initial vertex}
$\beg\,x$, and $\inv: D \to D$ is an involution which interchanges
every dart $x$ with its {\em inverse dart}, also denoted by $x^{-1}$.
The {\it neighbourhood} of a vertex $v$ is defined as the set
of darts that have $v$ for its initial vertex and the
{\it valence} of $v$ is the cardinality of the neighbourhood.
Note that a simple graph $\Gamma$ can be represented
as a pregraph by letting $V = \V(\Gamma)$ be the vertex-set of $\Gamma$,
letting $D = \{ (u,v) : u,v \in \V(\Gamma), u\sim_\Gamma v\}$ be
the arc-set of $\Gamma$, and by letting $\beg(u,v) = u$ and $\inv(u,v) = (v,u)$.
Notions such as morphism, isomorphism and automorphism of pregraphs are obvious
generalisations of those for (simple) graphs and precise definitions can be found in \cite{covers,elabcovers}.

The orbits of $\inv$ are called {\em edges}.
The edge containing a dart $x$ is called a {\em semi-edge} if $\inv\,x = x$,
a {\em loop} if $\inv\,x \neq x$ while
$\beg\,(x^{-1}) = \beg\,x$,
and  is called  a {\em link} otherwise.
The {\em endvertices of an edge} are the initial vertices of the darts contained in the edge.
Two links are {\em parallel} if they have the same endvertices.
When we present a pregraph as a drawing, the links are drawn in the usual way as a line
between the points representing its envertices, a loops is drawn as a closed curve at its unique endvertex and a semi-edge is drawn as a segment attached to its unique envertex.
The following lemma, which will serve as the starting point of our analysis of cubic tricirculants,
 is an easy exercise illustrating the definition of a pregraph.

\begin{figure}[h!]
\centering
\includegraphics[width=0.7\textwidth]{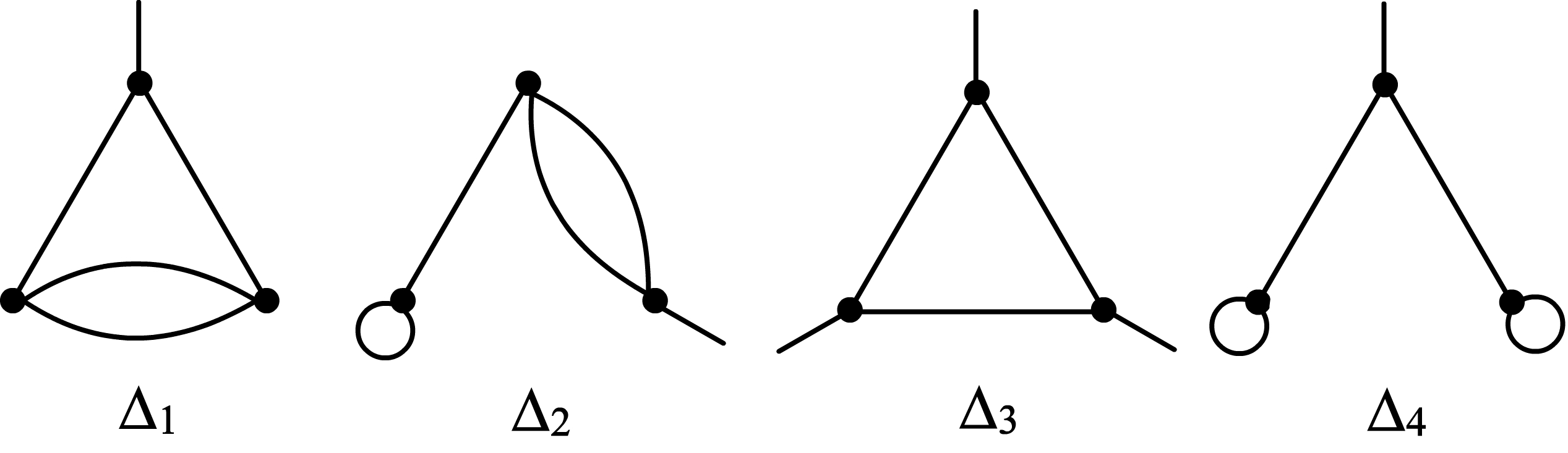}
\caption{Cubic pregraphs on three vertices without parallel semi-edges}
\label{fig:Delta}
 \end{figure}
 
\begin{lemma}
\label{lem:3D}
Up to isomorphism there are precisely four non-isomorphic cubic pregraphs on three vertices such that no vertex is the initial vertex of more than one semi-edge, namely
the pregraphs $\Delta_1, \Delta_2, \Delta_3$ and $\Delta_4$ depicted in Figure~\ref{fig:Delta}.
 \end{lemma}

Given a pregraph $\Gamma$ and a group $H\le \Aut(\Gamma)$ we define
the quotient $\Gamma/H$ to be the pregraph $(V',D'; \beg',\inv')$ where
$V'$ and $D'$ are the sets $\V(\Gamma)/H$ and $\D(\Gamma)/H$ of orbits of $H$ on $V$ and $D$, 
respectively, and $\beg'$ and $\inv'$ are defined by
$\beg' (x^H) = (\beg x)^H$ and $\inv' (x^H) = (\inv x)^H$ for every $x\in \D(\Gamma)$.
The mapping $\wp_H\colon \V(\Gamma) \cup \D(\Gamma) \to V' \cup D'$
that maps a vertex or a dart $x$ to its $H$-orbit $x^H$ is then an epimorphism of pregraphs, called {\it the quotient projection with respect to $H$}.
If $H$ acts semiregularly on $\V(\Gamma)$ (that is, if the vertex-stabiliser $H_v$ is trivial for every vertex $v\in \V(\Gamma)$), then the quotient projection $\wp_H$ is a {\em regular covering projection}  and in particular, it preserves the valence of the vertices.
Moreover, in this case the graph $\Gamma$ is isomorphic to the {\em derived covering graph} of 
$\Gamma/H$, a notion which we now define.

Let $\Gamma' = (V',D',\beg', \inv')$ be an arbitrary connected pregraph, let 
$N$ be a group and let $\zeta \colon D' \to N$ be a  mapping (called a 
{\em voltage assignment}) satisfying the condition 
$\zeta(x) = \zeta(\inv' x)^{-1}$ for every $x\in D'$. 
Then $\Cov(\Gamma',\zeta)$ is the graph with $D' \times N$ and $V' \times N$ 
as the sets of darts and vertices, respectively, and the functions 
$\beg$ and $\inv$ defined by $\beg(x,a) = (\beg' x,a)$ and 
$\inv(x,a) = (\inv' x, a \zeta(x))$. If, for a voltage assignment $\zeta\colon \D(\Gamma') \to N$,
there exists a a spanning tree $T$ in $\Gamma'$ such that $\zeta(x)$ is the trivial element of $N$
for every dart $x$ in $T$, then we say that $\zeta$ is {\it normalised};   
note that in this case $\Cov(\Gamma',\zeta)$ is connected if and only if the images of $\zeta$ generate $N$.
The following lemma is a well-known fact in the theory of covers:

\begin{lemma} 
\label{lem:normT}
Let $H$ be a group of automorphisms acting semiregularly on the vertex-set of a connected graph $\Gamma$, let $\Gamma'=\Gamma/H$ and let $T$ be a spanning tree in $\Gamma'$.
 Then there exists a voltage assignment $\zeta\colon \D(\Gamma') \to N$ which is normalised with respect to $T$ and such that $\Gamma \cong \Cov(\Gamma',\zeta)$.
\end{lemma}

\begin{figure}[h!]
\centering
\begin{tabular}{ccc}
\includegraphics[width=0.3\textwidth]{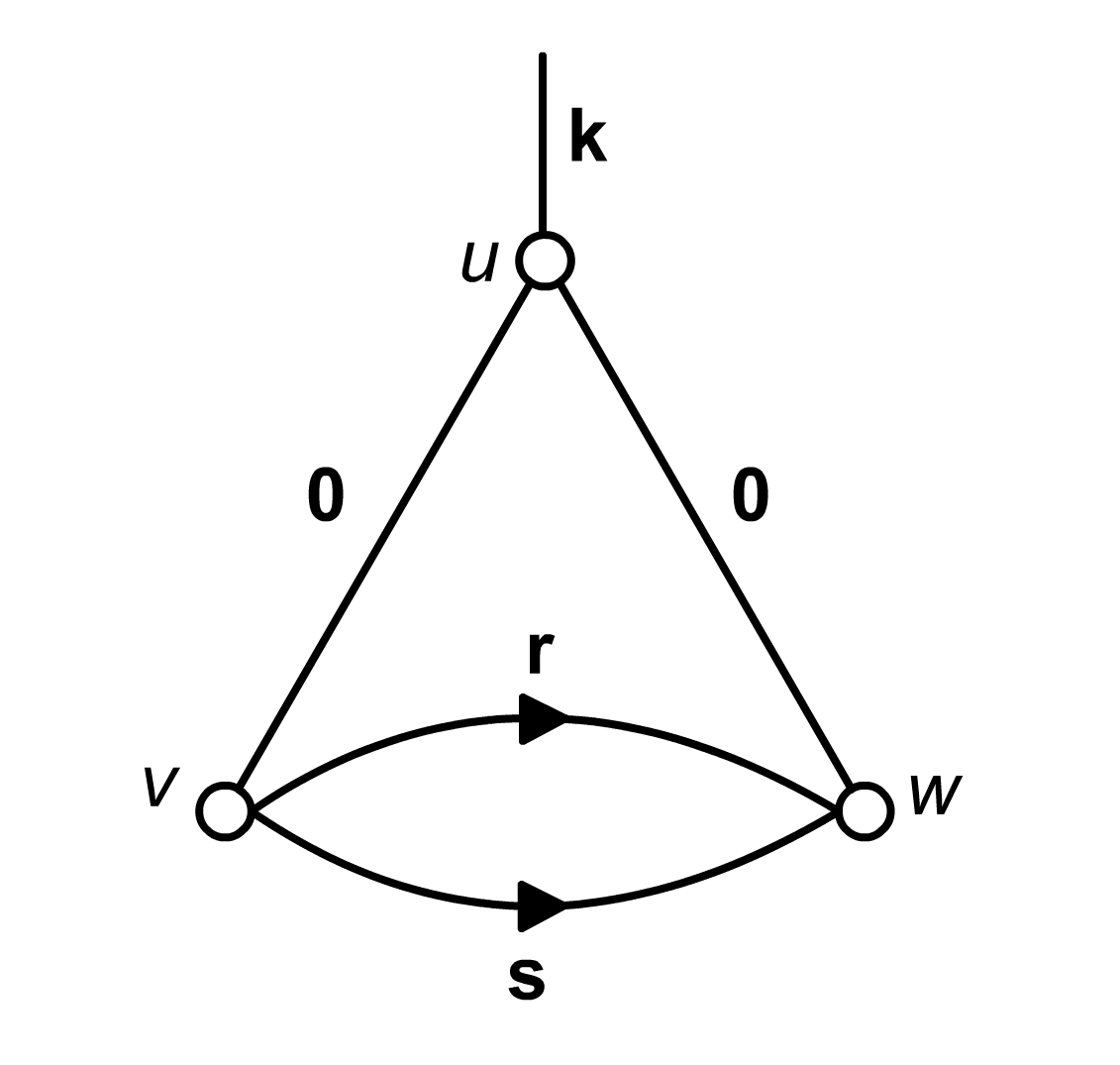} & &
\includegraphics[width=0.3\textwidth]{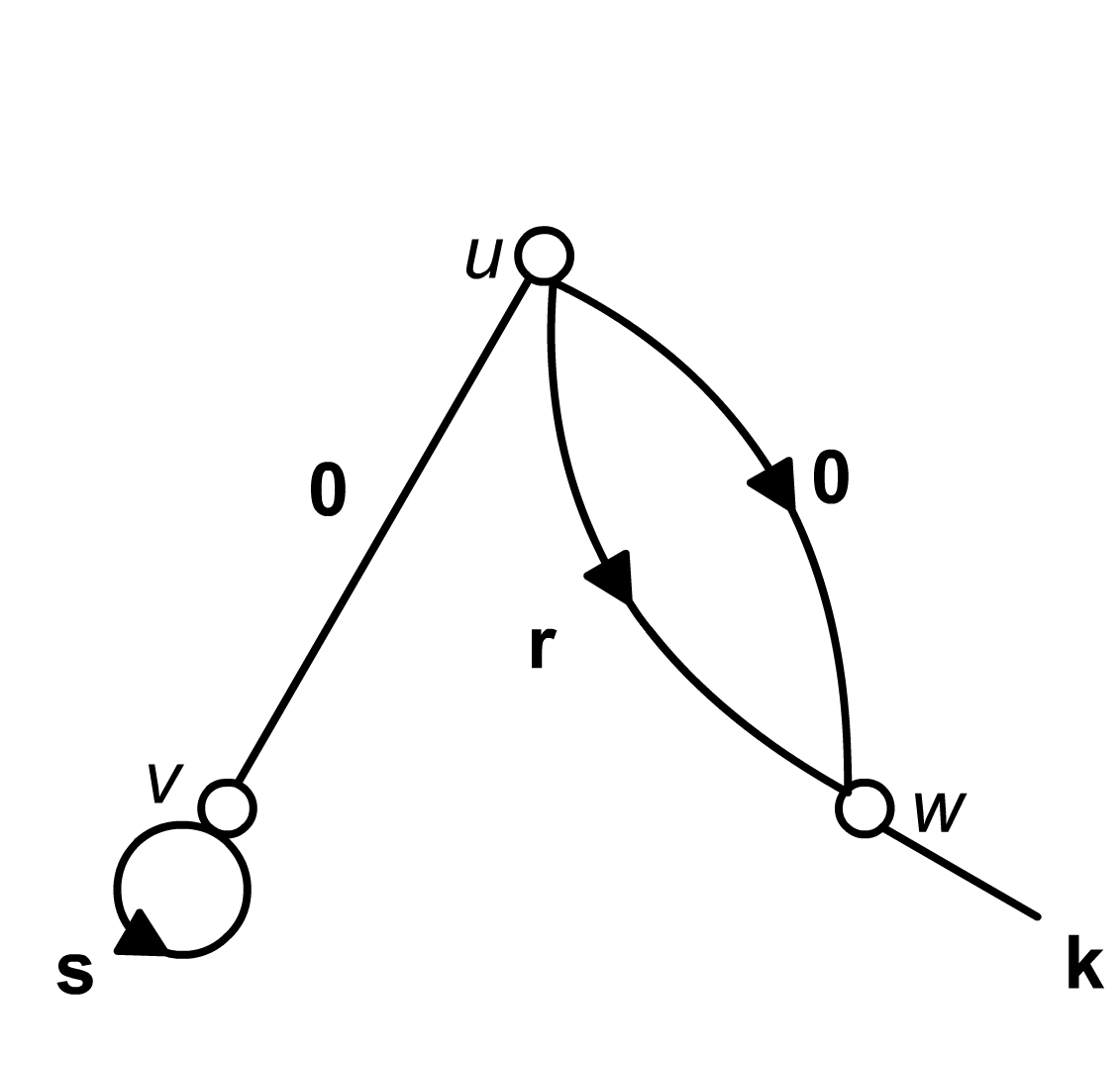} \\
$(\Delta_1,\zeta_1)$ & & $(\Delta_2,\zeta_2)$ \\  & \\
\includegraphics[width=0.3\textwidth]{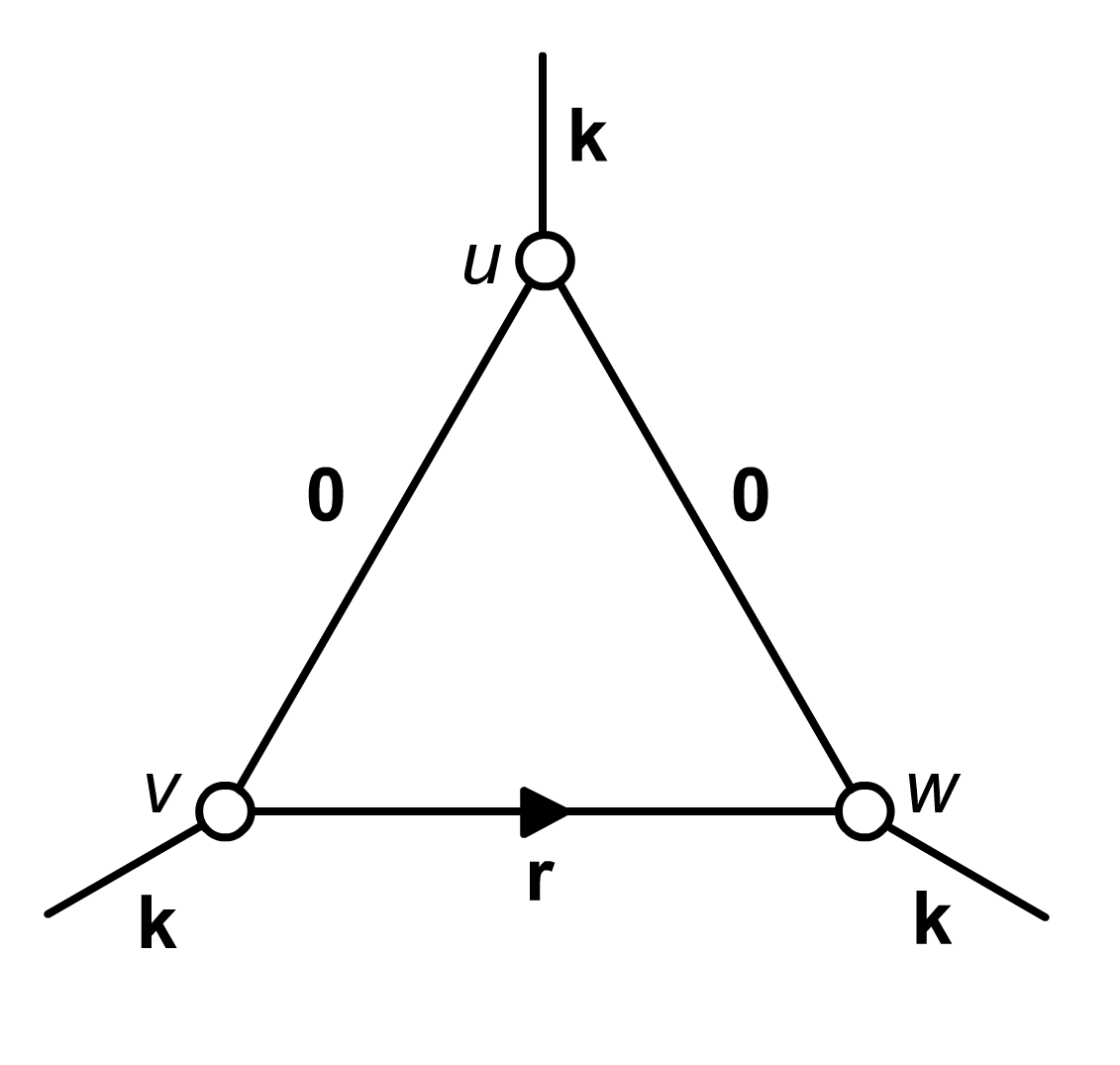} & &
\includegraphics[width=0.3\textwidth]{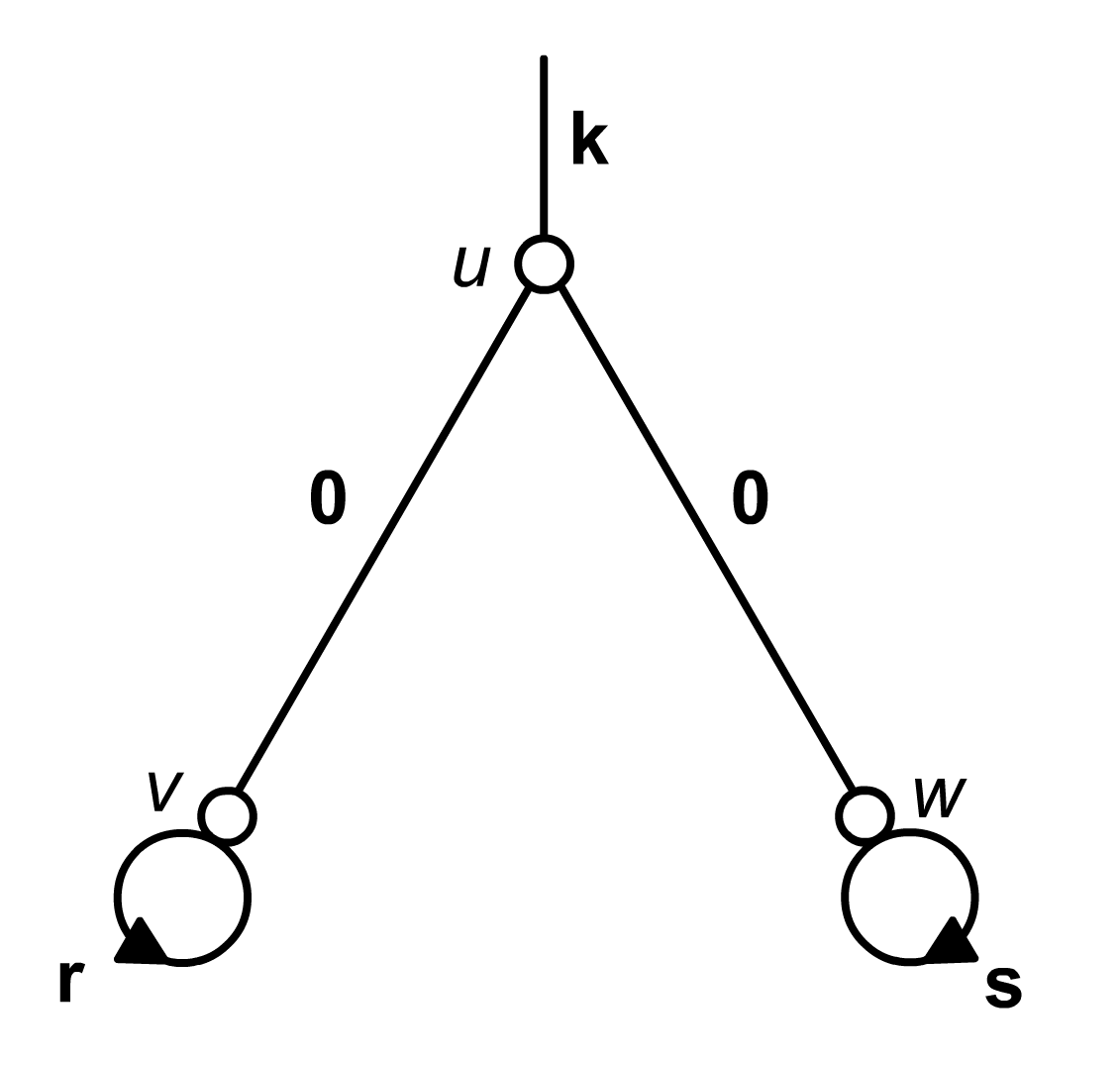} \\
$(\Delta_3,\zeta_3)$ & & $(\Delta_4,\zeta_4)$
\end{tabular}
\caption{The voltage assignments giving rise to cubic tricirculants.}
\label{fig:Ts}
 \end{figure}

Let us now move our attention back to cubic tricirculants. In the lemma
below, the voltage assignments $\zeta_i \colon \Delta_i \to \ZZ_n$ are given in Figure~\ref{fig:Ts}
by writing the value $\zeta_i(x)$ next to the drawing of each dart $x\in \D(\Delta_i)$.

\begin{lemma}
\label{lem:Ts}
Let $\Gamma$ be a cubic tricirculant, let $\rho$ be a corresponding tricirculant automorphism and
let $n$ be the order of $\rho$. Then $n=2k$ for some positive integer $k$ and there exist
elements $r,s \in \ZZ_n$ and $i\in \{1,2,3,4\}$ such that
 $\Gamma\cong \Cov(\Delta_i,\zeta_i)$ where $\zeta_i \colon \D(\Delta_i) \to \ZZ_n$ is
 the voltage assignment defined by Figure~\ref{fig:Ts}.
\end{lemma}

\begin{proof}
Note first that the quotient $\Gamma / \langle \rho \rangle$ is a pregraph with three vertices of valency $3$, one for each orbit under the action of $\langle \rho \rangle$.
Since $\rho$ is a semiregular automorphism of $\Gamma$ of order $n$, the group
$\langle \rho \rangle$ is isomorphic to $\ZZ_n$ and acts semiregularly on $\V(\Gamma)$. 
By Lemma~\ref{lem:normT},  $\Gamma \cong \Cov(\Gamma / \langle \rho \rangle,\zeta)$ for some voltage assignment $\zeta \colon \D(\Gamma / \langle \rho \rangle, \ZZ_n)$. Note that
by definition of voltage assignments, it follows that $\zeta(x) = -\zeta(x^{-1})$ for every $x \in \D(\Gamma / \langle \rho \rangle)$, implying that if $x$ is a semi-edge, then $\zeta(x)$ is an element of order at most $2$ in $\ZZ_n$, and since $\Gamma$ has no semi-edges, $\zeta(x)$ must in fact have order $2$ in $\ZZ_n$. Since every pregraph with three vertices in which every vertex has valence $3$ contains at least one semi-edge, it follows that $n=2k$ for some positive integer $k$, and moreover, $\zeta(x) = k$ for every semi-edge $x$ of $\Gamma / \langle \rho \rangle$. Since
$\Gamma$ has no parallel edges, this also implies that every vertex of $\Gamma / \langle \rho \rangle$ is the initial vertex of at most one semi-edge. By Lemma~\ref{lem:3D}, 
$\Gamma / \langle \rho \rangle \cong \Delta_i$ for some $i\in \{1,2,3,4\}$ and we may thus assume that $\Gamma\cong \Cov(\Delta_i,\zeta_i)$ for some $i\in \{1,2,3,4\}$ and some voltage assignment
$\zeta_i \colon \Delta_i \to \ZZ_n$. Finally, in view of Lemma~\ref{lem:normT}, we may assume
that $\zeta_i(x) = 0$ for every edge $x$ belonging to a chosen spanning tree of $\Delta_i$.
In particular, $\zeta_i$ can be chosen  as shown in Figure~\ref{fig:Ts}.
\end{proof}

A cubic tricirculant isomorphic to $\Cov(\Delta_i,\zeta_i)$ is said to be of {\em Type $i$}.
Note that in principle, a cubic vertex-transitive tricirculant could be of more than one type;
however, it can be shown that this is not the case except for the complete bipartite graph
$K_{3,3}$ which is both of Type 1 as well as of Type 3.

The following sections are devoted to the analysis of the tricirculants graphs arising from the voltage assignments $\zeta_i$, and in particular, to determining sufficient and necessary conditions for vertex-transitivity. We will consider the graphs of order at most $48$ separately in Section~\ref{sec:small}
and as for the graphs of larger order, we will show that no graph of Type 4 is vertex-transitive, that
Type 3 yields prisms and M\"obius ladders, and that Types 1 and 2 each yield one infinite family
of vertex-transitive graphs, namely the graphs $\textrm{X}(k)$ defined in Definition~\ref{def:X} 
and the graphs $\textrm{Y}(k)$ defined in Definition~\ref{def:Y}.

Finally, to facilitate discussion in the forthcoming sections, we specify the notion of a walk.
A walk $\omega$ in a (pre)graph is a sequence of darts $(x_1,x_2,...,x_n)$ such that $\beg(x_{i+1}) = \beg(x_i^{-1})$ for all $i \in \{1..,n-1\}$; that is, the initial vertex of $x_{i+1}$ is the end vertex of $x_i$. We say that $\omega$ is closed if the initial vertex of $x_1$ is the end vertex of $x_n$. If additionally $\beg(x_i) \neq \beg(x_j)$ for all $i \neq j$, then we say that $\omega$ is a cycle. We define the inverse of $\omega$ as the walk $\omega^{-1}=(x_n^{-1},x_{n-1}^{-1},\ldots, x_1^{-1})$.
We say that a walk $\omega$ is a {\em reduced walk} if no two consecutive arcs are inverse to one another (if $\omega$ is closed, we will consider the first and last arcs of $\omega$ to be consecutive). It is pertinent to point out that, from this definition, all walks and cycles are directed and have an initial vertex. Therefore, no reduced walk is equal to its inverse.

\section{Graphs of small order}
\label{sec:small}

To make this classification as general and as neat as possible, we will treat cubic tricirculant graphs of small order separately, as if we allow graphs that are ``too small'', special cases and exceptions will inevitably occur. Since a cubic tricirculant must have order $6k$, for some positive $k$, we present, in Table \ref{table:small}, the complete list of all cubic vertex-transitive tricirculants of order at most $48$, along with their Types, obtained from the census \cite{census} of cubic vertex-transitive graphs. Each graph in the census is labelled with a pair of integers in square brackets, the first being the order of the graph and the second being just an identifier. \\

\begin{table}[H]
\begin{tabular}{|c c| c c| c c| c c| c c|}
\hline
 Label & Type & Label & Type & Label & Type & Label & Type & Label & Type \\
\hline
 [6,1] & $1,3$ & [18,1] & $3$ & [24,1] & $3$ & [30,7] & $3$ & [42,4] & $1$ \\
\hline
 [6,2] & $3$ & [18,3] & $3$ & [30,1] & $2$ & [30,8] & $4$ & [42,9] & $3$\\
\hline
 [12,1] & $3$ & [18,4] & $2$ & [30,2]  & $1$ & [36,1] & $3$ & [42,10] & $3$ \\
\hline
 [12,2] & $1$ & [18,5] & $1$ & [30,5]  & $3$ & [42,3] & $2$ & [48,1] & $3$\\
\hline
\end{tabular}
\caption{Graphs of small order}
\label{table:small}
\end{table}
The graphs [6,1], [18,4] and [30,8] are arc-transitive and correspond to $K_{3,3}$, $F018A$ (the Pappus graph) and $F030A$ (Tutte's $8$-cage) respectively (see \cite{foster,symtric} for further details). The fourth arc-transitive cubic tricirculant, denoted $F054A$ in \cite{foster}, 
does not appear in Table \ref{table:small}, as it has $54$ vertices; it corresponds to the graph $\textrm{Y}(9)$  defined in Section~\ref{Type2}. Since a cubic vertex-transitive graph that is also edge-transitive is automatically arc-transitive \cite{tutte}, we have that no vertex-transitive cubic tricirculant is edge-transitive, with the exception of the four aforementioned graphs.

For the remainder of this paper, all cubic tricirculants will have order at least $54$.

\section{Type 1} \label{type1}

\begin{figure}[h!]
\centering
\includegraphics[width=0.3\textwidth]{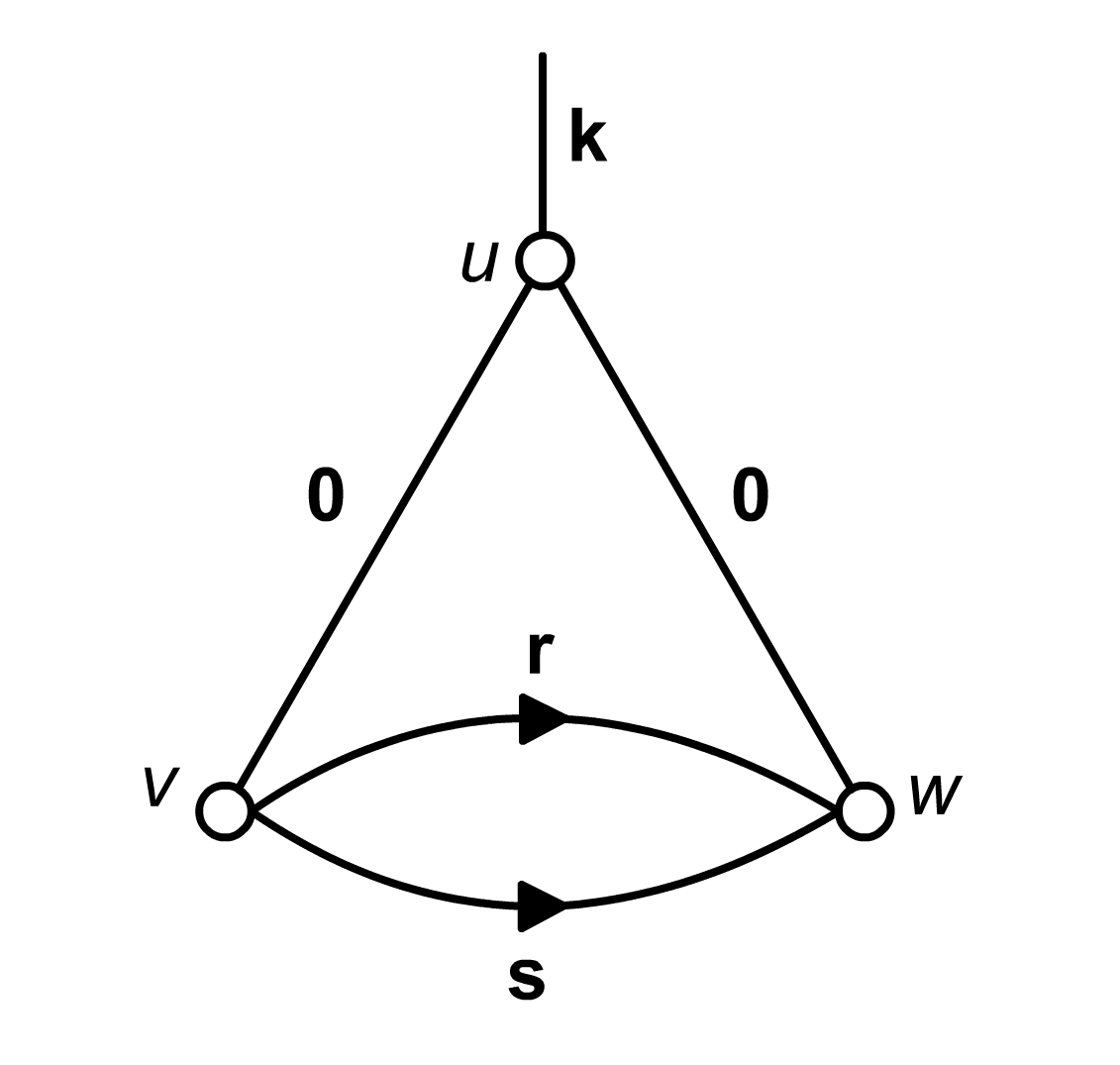}
\caption{The voltage assignment giving rise to the graph $T_1(k,r,s)$.}
\label{fig:Delta1}
 \end{figure}

Let $k \geq 9$ be an integer, and let $r$ and $s$ be two distinct elements of $\ZZ_{2k}$. Let $T_1(k,r,s)$ be the covering graph arising from $\zeta_1$ given in Figure~\ref{fig:Ts}; for convenience, we repeat the drawing here (see  Figure~\ref{fig:Delta1}).

By the definition of a derived cover, the vertices of 
$T_1(k,r,s)$ are then the pairs $(x,i)$ where $x \in \{u,v,w\}$
is a vertex of $\Delta_1$ and $i\in \ZZ_{2k}$. Note that $T_1(k,r,s)$ is connected if and only if $k$, $r$ and $s$ generate $\ZZ_{2k}$, that is, if and only if $\gcd(k,r,s)=1$.

To simplify notation, we write
$u_i$ instead of $(u,i)$ and similarly for $v_i$ and $w_i$. 
Then $U=\{u_0,u_1,...,u_{2k-1}\}$, $V=\{v_0,v_1,...,v_{2k-1}\}$ and $W=\{w_0,w_1,...,w_{2k-1}\}$ are the respective fibres of vertices $u$, $v$ and $w$,  and the edge-set of $T_1(k,r,s)$
is then the union of the sets
\begin{eqnarray}
E_K &=& \{u_iu_{i+k}: i \in \ZZ_{2k}\}, \cr
E_0 &=& \{u_iv_i: i \in \ZZ_{2k}\} \cup \{u_iw_i: i \in \ZZ_{2k}\}, \cr
E_R &=& \{v_iw_{i+r}: i \in \ZZ_{2k}\}, \cr 
E_S &=& \{v_iw_{i+s}: i \in \ZZ_{2k}\}. \cr
\end{eqnarray}

\begin{definition}
\label{def:X}
For an odd positive integer $k$, let
$$
r^* = \left\{ \begin{array}{ll} \frac{k+3}{2}, & \hbox{ if }  k \equiv 1\> (\hbox{\rm mod}\> 4) \\
   \frac{k+3}{2} + k, & \hbox{ if }  k \equiv 3\> (\hbox{\rm mod}\> 4) \end{array}\right.
$$
and let $\textup{X} (k) = T_1(k,r^*,1)$.
\end{definition}

We can now state the main theorem of this section. 

\begin{theorem}
\label{the:t1}
Let $\Gamma$ be a connected cubic tricirculant of Type 1 with at least $54$ vertices.
Then $\Gamma$ is vertex-transitive if and only if it is isomorphic to $\textup{X}(k)$ for some odd integer $k$, $k\ge 9$.
In this case, $\Gamma$ is a bicirculant if $3 \nmid k$ and is isomorphic to the Generalized Petersen graph $\textrm{GP}(3k,k+(-1)^\alpha)$ where $\alpha \in \{1,2\}$ and $\alpha \equiv k$ (\mod$3$), and thus a bicirculant.
\end{theorem}

The rest of the section is devoted to the proof of Theorem~\ref{the:t1}.
Clearly, as a tricirculant of Type 1, $\Gamma$ is isomorphic to $T_1(k,r,s)$ for some integer $k$
and elements $r,s\in \ZZ_{2k}$ (see Lemma~\ref{lem:Ts}). We shall henceforth assume that
$\Gamma = T_1(k,r,s)$ for some $k\ge 9$ and that $\Gamma$ is connected.
For a symbol $X$ from the set of symbols $\{0,R,S,K\}$, edges in $E_X$ will be called edges of type $X$, or simply $X$-edges. 

Define $\rho$ as the permutation given by $\rho(u_i)=u_{i+1}$, $\rho(v_i)=v_{i+1}$  and $\rho(w_i)=w_{i+1}$, and note that $\rho$ is a $3$-circulant automorphism of $T_1(k,r,s)$. 
%
Further, observe that 
\begin{equation}
\label{eq:rs}
T_1(k,r,s) \cong  T_1(k,s,r),
\end{equation}
and if $a \in \ZZ$ is such that $\gcd(2k,a)=1$, then
\begin{equation}
\label{eq:auto}
T_1(k,ar,as) \cong T_1(k,r,s).
\end{equation}

\begin{lemma}
\label{lem: xvt}
Let $k$ be an odd integer, $k\ge 9$, and let $r^*$ and $\textup{X}(k)$ be as in Definition~\ref{def:X}.
Then the graph $X(k)$  is vertex-transitive. 
\end{lemma}

\begin{proof}
Recall that $\textup{X}(k) \cong T_1(k,r^*,1)$. Define $\phi$ as the mapping given by:

\begin{center}
\begin{tabular}{llll}
$u_i \mapsto w_{i-r^*+2}$, &  $w_i \mapsto v_{i-2r^*+2}$, & $v_i \mapsto u_{i-r^*+2}$  & if $i$ is even; \\
$u_i \mapsto v_{i+r^*-2}$,  &  $v_i \mapsto w_{i+2r^*-2}$, & $w_i \mapsto u_{i+r^*-2}$ & if $i$ is odd.
\end{tabular}
\end{center} 

That $\phi$ is indeed an automorphism follows from the fact that $k$ is odd, $r^*$ is even and the congruence $k+3-2r^* \equiv 0$ holds. Since $\phi$ transitively permutes the three $\langle \rho\rangle$-orbits $U$, $V$ and $W$, the group $\langle \rho, \phi\rangle$ acts transitively on the vertices of $T_1(k,r^*,1)$.
\end{proof}

Now, recall that $\Gamma$ is connected and equals $T_1(k,r,s)$, for some $k\ge 9$, $r,s \in \ZZ_{2k}$. We may also assume that $\Gamma$ is vertex-transitive, however, for some of the results in this section vertex-transitivity is not needed. When possible, we will not assume the graph to be vertex-transitive, but rather only to have some weaker form of symmetry, that we will define in the following paragraphs.

A simple graph $\Gamma$ is said to be {\em $c$-vertex-regular}, for some $c \geq 3$, if there are the same number of $c$-cycles through every vertex of $\Gamma$.

For an edge $e$ of $\Gamma$ and a positive integer $c$ denote by $\epsilon_c(e)$ the number of $c$-cycles that pass through $e$. For a vertex $v$ of $\Gamma$, let $\{e_1,e_2,e_3\}$ be the set of edges incident to $v$ ordered in such a way that $\epsilon_c(e_1) \leq \epsilon_c(e_2) \leq \epsilon_c(e_3)$. The triplet $(\epsilon_c(e_1),\epsilon_c(e_2),\epsilon_c(e_3))$ is then called the {\em $c$-signature} of $v$. If for a $c \geq 3$ all vertices in $\Gamma$ have the same $c$-signature, then we say $\Gamma$ is {\em $c$-cycle-regular} and the signature of $\Gamma$ is the signature of any of its vertices. Note for every $c\ge 3$,
a vertex-transitive graph is necessarily $c$-cycle-regular, and every $c$-cycle-regular graph is $c$-vertex-regular.

If $c$ equals the girth of the graph, then following 
\cite{PotVid}, a 
$c$-cycle-regular graph will be called 
{\em girth-regular}.


\begin{lemma}
\label{lem_rsnok}
If $\Gamma$ is $4$-vertex-regular, then neither $r$ nor $s$ equals $k$.
\end{lemma}

\begin{proof}
Suppose $r=k$. Since $\Gamma$ has no parallel edges, we see that $s \neq k$. Observe that 
$(u_0,v_0,w_r,u_r)$ and $(u_0,w_0,v_{-r},u_{-r})$ are the only $4$ cycles of $\Gamma$ through $u_0$. 

 Meanwhile, there exists a unique $4$-cycle through $v_0$, namely $(v_0,w_r,u_r,u_0)$, which contradicts $\Gamma$ being $4$-vertex-regular. 
Hence $r\not = k$, and view of the isomorphism $T_1(k,s,r) \cong T_1(k,r,s)$ 
 (see (\ref{eq:rs})), this also shows that $s\not = k$.
\end{proof}
 
 For the sake of simplicity denote a dart in $\Delta_1$ starting at vertex $a$, pointing to vertex $b$ and having voltage $x$ by $(ab)_x$. Recall that a walk in $\Delta_1$ is a sequence of darts $(x_1,x_2,...,x_n)$, for instance $((vw)_r,(wu)_0,(uu)_k)$ is a walk of length $3$. However, since $\Gamma$ is a simple graph a walk in $\Gamma$ will be denoted as a sequence of vertices, as it is normally done.

\begin{lemma}
\label{lem:rsno0}
If $\Gamma$ is $8$-cycle-regular, then neither $r$ nor $s$ equals $0$.
\end{lemma}

\begin{proof}
Suppose $r=0$. Since $\Gamma$ has no parallel edges, we see that $s \neq 0$. Suppose there is an $8$-cycle $C$ in $\Gamma$. Such a cycle, when projected to $\Delta_1$, yields a reduced closed walk $\omega$ in $\Delta_1$ whose $\zeta_1$-voltage is $0$ . Note that $\omega$ cannot trace three darts with voltage $0$ consecutively, as this would lift into a $3$-cycle contained in $C$. This implies $\omega$ necessarily visits the dart $(vw)_s$ or its inverse at least once. Furthermore, since $\gcd(k,s)=1$ and $8<9 \leq k$, $\omega$ must trace $(wv)_{-s}$ as many times as it does $(vw)_{s}$. By observing Figure~\ref{fig:Delta1}, the reader can see that if $\omega$ traces the dart $(vw)_s$, then it must also trace the semi-edge $(uu)_k$ before tracing $(wv)_{-s}$, as $\omega$ cannot trace $3$ consecutive darts with voltage $0$. Since $\omega$ has net voltage $0$, it necessarily traces $(uu)_k$ an even amount of times. Moreover, $\omega$ must trace a dart in $\{(uv)_0,(vu)_0,(uw)_0,(wu)_0\}$ immediately before and immediately after tracing $(uu)_k$. Hence, $\omega$ must visit the set $\{(vw)_s,(wv)_{-s}\}$ at least twice; the semi-edge $(uu)_k$ at least twice; and the set $\{(uv)_0,(vu)_0,(uw)_0,(wu)_0\}$ at least $4$ times. This already amounts to $8$ darts, none of which is $(vw)_r$ or its inverse. Therefore, no $8$-cycle in $\Gamma$ visits an $R$-edge.
However, for $X \neq R$ there is at least one $8$-cycle through every $X$-edge, as $(u_0, v_0, w_s, u_s, u_{s+k}, w_{s+k}, v_k, u_k)$ and $(u_0, w_0, v_{-s}, u_{-s}, u_{k-s}, v_{k-s}, w_k, u_k)$ are $8$-cycles in $\Gamma$. This contradicts our hypothesis of $\Gamma$ being $8$-cycle-regular. The proof when $s=0$ is analogous.
\end{proof}

Now, observe that $(u_0,u_k,v_k,w_{k+r},u_{k+r},u_r,w_r,v_0,u_0)$ is a cycle of length $8$ in $\Gamma$ starting at $u_0$. In what follows we will study the $8$-cycle structure of $\Gamma$ to determine conditions for vertex-transitivity. Recall that each $8$-cycle in $\Gamma$ quotients down into a closed walk of length $8$ having net voltage $0$ in $\Delta_1$. We can thus, provided we are careful, determine how many $8$-cycles pass through any given vertex of $\Gamma$ by counting closed walks of length $8$ and net voltage $0$  in the quotient $\Delta_1$. Note that only closed walks that are reduced will lift into cycles of $\Gamma$. Therefore, we may safely ignore non-reduced walks and focus our attention exclusively on reduced ones. Define $\mathcal{W}_8$ as the set of all reduced closed walks of length $8$ in $\Delta_1$. 

Every element of $\mathcal{W}_8$ having net voltage $0$ lifts into a closed walk of length $8$ in $\Gamma$. In principle, the latter  walk might not be a cycle of $\Gamma$. However, we will show that this never happens in our particular case.
For this, it suffices to show that $\Gamma$ does not admit cycles of length $4$ or smaller, as any closed walk of length $8$ that is not a cycle is a union of smaller cycles, one of which will have length smaller or equal to $4$.

\begin{lemma}
\label{lem:girth5}
$\Gamma$ has girth at least 5.
\end{lemma}

\begin{proof}
Suppose $\Gamma$ has a $3$-cycle. Since $\Gamma$ is vertex-transitive, this would mean that there is a reduced closed walk with net voltage $0$ of length $3$ in $\Delta_1$ that visits $u$. From Figure \ref{fig:Delta1} we get that such a closed walk must be either $((uv)_0,(vw)_r,(wu)_0)$, $((uv)_0,(vw)_s,(wu)_0)$ or one of their inverses. This would imply that either $r=0$ or $s=0$, a contradiction. Similarly, the existence of $4$-cycles in $\Gamma$ would imply there is a closed walk with net voltage $0$ of length $4$ visiting $u$. The only such walks are $((uv)_0,(vw)_r,(wu)_0,(uu)_k))$, $((uu)_k,(uv)_0,(vw)_r,(wu)_0)$, $((uv)_0,(vw)_s,(wu)_0,(uu)_k))$, $((uv)_0,(vw)_s,(wu)_0,(uu)_k))$ and their inverses. In all eight cases, we have that $r=k$ or $s=k$, again a contradiction.
\end{proof}

Now, let $N$ be the set of the net voltages of walks in $\mW_8$, expressed as linear combinations of $k$, $r$ and $s$, where 
we view $r$ and $s$ as indeterminants over $\ZZ_{2k}$. For instance, the closed walk $((uv)_0,(vw)_r,(wv)_{-s},(vw)_r,(wu)_0,(uv)_0,(vw)_r,(wu)_0)$ has net voltage $3r - s$.
For  an element $\nu \in N$, let $\mW(\nu)$ be the set of walks in $\mW$ with net voltage $\nu$ and for  a
vertex $x$ of $\Delta_1$, let $\mW_8^x(\nu)$ be the set of walks in $\mW(\nu)$ starting at $x$.

 For each $\nu\in N$, we have computed the number of elements in $\mW_8^u(\nu)$ and in $\mW_8^v(\nu)$. 
  The result is displayed in Table~\ref{table:1}. Notice that, if $\omega \in \mW_8(\nu)$, then $\omega^{-1}\in \mW_8^x(-\nu)$.
   This means $\omega$ has voltage $0$ if and only if $\omega^{-1}$ does too. For this reason, and since we are interested in walks with net voltage $0$, we have grouped walks with net voltage $\nu$ along with their inverses having net voltage $-\nu$.  This computation is straightforward but somewhat lengthy. For this reason and to avoid human error, we have done it with the help of a computer programme.

\begin{table}
\begin{tabular}{|c| c| c| c|}
\hline
 Label & Net voltage & Starting at $u$ & Starting at $v$   \\
 \hline
 $\textrm{I}$ & $0$ & $12$ & $10$ \\
 $\textrm{II}$ & $\pm(2r)$ & $8$ & $8$  \\
 $\textrm{III}$ & $\pm(2s)$ & $8$ & $8$  \\
 $\textrm{IV}$ & $\pm(r+s)$ & $8$ & $4$  \\
 $\textrm{V}$ & $\pm(r-s)$ & $8$ & $4$ \\
 $\textrm{VI}$ & $\pm(k+2r-s)$ & $12$ & $10$  \\
 $\textrm{VII}$ & $\pm(k+2s-r)$ & $12$ & $10$  \\
 $\textrm{VIII}$ & $\pm(k+3r-2s)$ & $4$ & $6$  \\
 $\textrm{IX}$ & $\pm(k+3s-2r)$ & $4$ & $6$  \\
 $\textrm{X}$ & $\pm(3r-s)$ & $4$ & $6$  \\
 $\textrm{XI}$ & $\pm(3s-r)$ & $4$ & $6$  \\  
 $\textrm{XII}$ & $\pm(2r-2s)$ & $4$ & $6$  \\
 $\textrm{XIII}$ & $\pm(4r-4s)$ &$0$ & $2$\\
 \hline
\end{tabular}
\caption{Net voltages of closed walks of length $8$ in $\Delta_1$.}
\label{table:1}
\end{table}


From the first row of Table \ref{table:1} we know that there are always $12$ walks starting at $u$ in $\mW_8$ that have net voltage $0$, regardless of the values of $r$ and $s$. Similarly, there are always $10$ such walks starting at $v$. Since $\Gamma$ is vertex-transitive, the amount of walks in $\mW_8$ having net voltage $0$ (after evaluating $r$ and $s$ in $\ZZ_{2k}$) starting at $u$ and those starting at $v$ must be the same. It follows that, for some $\nu \in N$ with $|\mW_8^v(\nu)| > |\mW_8^u(\nu)|$, the equation $\nu \equiv 0$ (\mod $2k$) holds when we evaluate $r$ and $s$ in $\ZZ_{2k}$. We thus see that at least one expression in $\textrm{VIII}$--$\textrm{XIII}$ is congruent to $0$ modulo $2k$. In fact, if $\Gamma$ is vertex-transitive, then at most one of these expressions can be congruent to $0$.

\begin{lemma}
\label{lem:eqtype1}
If $\Gamma$ is $4$-vertex-regular and $8$-cycle-regular, then exactly one of the following equations modulo $2k$ holds:
\begin{eqnarray}
\label{eq1} 3s-2r+k \equiv 0\\
\label{eq2} 3r-2s+k \equiv 0\\
\label{eq3} 3r-s \equiv 0\\
\label{eq4} 3s-r \equiv 0\\
\label{eq6} 4r-4s \equiv 0
\end{eqnarray}
\end{lemma}

\begin{proof}
We will slightly abuse notation and, for a label $x \in \{\textrm{I},\textrm{II},...,\textrm{XIII}\}$, we will refer by $x$ to the congruence equation modulo $2k$ obtained by making the expression labelled $x$ in Table \ref{table:1} congruent to $0$. For instance, the equation $2r \equiv 0$ will be referred to as $\textrm{II}$. 

First note that if $\textrm{II}$ or $\textrm{III}$ holds then $r=k$ or $r=k$, contradicting Lemma \ref{lem_rsnok}. If $\textrm{V}$ holds then $(v_0,w_r,v_0)$ is a $2$-cycle in $\Gamma$, which is not possible. Similarly, $\textrm{XII}$ implies the existence of a $4$-cycle in $\Gamma$, contradicting Lemma \ref{lem:girth5}.  Now, if $\textrm{XIII}$ holds, then either $2r-2s \equiv 0$, which is not possible, or $2r-2s+k \equiv 0$.
 
  If $\textrm{VI}$ holds, then neither $\textrm{X}$ nor $\textrm{XI}$ can hold as this would imply $2r-2s \equiv 0$ (subtracting $\textrm{IV}$ from $\textrm{X}$ or from $\textrm{XI}$, respectively). Therefore, $\textrm{IV}$ excludes $\textrm{X}$ and $\textrm{XI}$. If $\textrm{IV}$ and $\textrm{I}$ are the only equations to hold, we would have $6$ more elements of $\mW_8$ through $u$ than through $v$. This means that if $\textrm{IV}$ holds then necessarily $\textrm{VII}$, $\textrm{IX}$ and $\textrm{XII}$ also hold. However, $\textrm{XII}$ implies $2r-2s+k \equiv 0$ and subtracting this from $\textrm{VIII}$ we get $r=0$, a contradiction. Hence, $\textrm{VII}$, $\textrm{IX}$ and $\textrm{XII}$ cannot all hold at the same time and so $\textrm{IV}$ can never hold. 
  
Suppose $\textrm{VI}$ holds, then one of the equations in $\{\textrm{VIII},\textrm{IX},\textrm{X}, \textrm{XI}, \textrm{XIII}\}$ must also hold. Note that $\textrm{VI}$ and $\textrm{VIII}$ imply $\textrm{III}$; $\textrm{VI}$ and $\textrm{IX}$ imply $\textrm{V}$; $\textrm{VI}$ and $\textrm{X}$ imply $r=k$; $\textrm{VI}$ and $\textrm{XIII}$ imply $s=k$. Therefore only $\textrm{XI}$ can hold. But $\textrm{VI}$ and $\textrm{XI}$ imply $\textrm{VII}$, and so we would have $40$ elements of $\mW_8$ starting at $u$ but only $36$ starting at $v$. This would contradict vertex-transitivity. Thus $\textrm{VI}$ cannot hold. An analogous reasoning, where the roles of $r$ and $s$ are interchanged, shows that $\textrm{VII}$ cannot hold. 
We have shown that the only equations that can hold are in $\{\textrm{VIII},\textrm{IX},\textrm{X}, \textrm{XI}, \textrm{XIII}\}$. However, if two or more of these equations hold, then we would have more elements of $\mW_8$ through $v$ than we would have through $u$. It follows that exactly one equation in $\{\textrm{VIII},\textrm{IX},\textrm{X}, \textrm{XI}, \textrm{XIII}\}$ holds.
\end{proof}

Lemma \ref{lem:eqtype1} tells us, for each of the $5$ possible equations, exactly which walks in $\mW_8$ will lift to $8$-cycles and so we can count exactly how many $8$-cycles pass through a given vertex or edge of $\Gamma$. For instance, if $3s-2r+k \equiv 0$, then there are $16$ walks in $\mW_8$ through $u$ that will lift to an $8$-cycle (see Table \ref{table:1}). Since a closed walk in $\mW_8$ and its inverse lift to the same cycle in $\Gamma$, we see that there are $8$ cycles of length $8$ through every vertex in $\Gamma$. Moreover, there are exactly $5$ cycles of length $8$ through every edge of type $R$ or $0$ while there are $6$ such cycles through each edge of type $K$ or $S$. It follows that the $8$-signature of a vertex $x$ of $\Gamma$ is $(5,5,6)$ whenever $3s-2r+k \equiv 0$ and thus, in this case, $\Gamma$ is $8$-cycle-regular with $8$-signature $(5,5,6)$. Table \ref{table:2} shows the number of $8$-cycles through each vertex and through each edge, depending on its type, for each of the $5$ cases described in Lemma \ref{lem:eqtype1}.

\begin{table}
\begin{tabular}{|c| c| c| c| c| c|}
\hline
 Congruence & $0$-edge & $R$-edge & $S$-edge & $K$-edge & $8$-signature    \\
\hline
 $3s -2r +k$ & $5$ & $5$ & $6$ & $6$ & $(5,5,6)$ \\
 $3r -rs +k$ & $5$ & $6$ & $5$ & $6$ & $(5,5,6)$\\
 $3r - s$ & $6$ & $6$ & $4$ & $4$ & $(4,6,6)$\\
 $3s -r$ & $6$ & $4$ & $6$ & $4$ &  $(4,6,6)$\\
 $4r-4s$ &$4$ & $4$ & $4$ & $4$ & $(4,4,4)$\\
 \hline
 \end{tabular}
\caption{Number of $8$-cycles through each edge-type of $T_1(k,r,s)$}
\label{table:2}
\end{table}

We will now show that equations \ref{eq3}, \ref{eq4} and \ref{eq6} cannot hold when $\Gamma$ is vertex-transitive. This will be proved in Lemmas~\ref{lem:4.6}, \ref{lem:4.7} and \ref{lem:4.8}. But first we need to show a result about the subgraph induced by the $0$- and $R$-edges.

Recall that that $\Gamma=T_1(k,r,s)$ with $k\ge 9$ and that it is connected.
Henceforth we also assume that $\Gamma$ is vertex-transitive.
Denote by $\Gamma_{0,R}$ the subgraph that results from deleting the edges of type $S$ and $K$ from $\Gamma$. Equivalently, $\Gamma_{0,R}$ is the subgraph induced by $0$- and $R$-edges. We see that $\Gamma_{0,R}$ is $2$-valent and that it has $\gcd(2k,r)$ connected components, each of which is a cycle of length $6k/\gcd(2k,r)$. 

Since $\Gamma_{0,R}$ only has edges of type $0$ and $R$, any two vertices in the same connected components must have indices that differ in a multiple of $r$. It is not hard to see that, indeed, each connected component consist precisely of all those vertices whose indices are congruent modulo $\gcd(2k,r)$. We now prove an auxiliary result that will used both in the case where (\ref{eq3}) as well as in the case where
(\ref{eq1}) or  (\ref{eq2}) holds. 

Suppose that no automorphism of $\Gamma$ maps an $S$-edge to a $0$-edge.
Since $\Gamma$ is vertex-transitive, there is an automorphism that maps
a vertex from $V$ to a vertex in $U$. Such an automorphism then maps an $S$-edge to a $K$-edge. Since all $S$-edges as well as all the $K$-edges are in the same
orbit, this then implies that $E_S \cup E_K$ forms a single edge-orbit of $\Aut(\Gamma)$.

\begin{lemma}
\label{lem_blocks}
Suppose that no automorphism of $\Gamma$ maps an $S$-edge to a $0$-edge, then $\Gamma_{0,R}$ is disconnected, and the set of connected components induce a block system for the vertex-set of $\Gamma$.
\end{lemma}

\begin{proof}
By the discussion above the lemma, $E_S \cup E_K$ forms a single edge-orbit in $\Gamma$.
Now note that $\Aut(\Gamma)$ also acts as a group of automorphism on $\Gamma_{0,R}$ and that this action is transitive, since the edges removed consist of an edge-orbit of $\Gamma$. It follows that connected components of $\Gamma_{0,R}$ form a block system for the vertex set of $\Gamma$. Now, suppose $\Gamma_{0,R}$ consists of a single cycle of length $6k$, $(w_0,u_0,v_0,w_r,u_r,v_r \dots v_{(2k-1)r})$. Notice that the vertex antipodal to $u_0$ in this cycle is $u_{kr}$ which is the same as $u_k$, as $r$ is odd. We see that edges of type $K$ in $\Gamma$ join vertices that are antipodal in $\Gamma_{0,R}$, while $S$-edges do not. Since this $6k$-cycle is a block of imprimitivity, $K$-edges can only be mapped into $K$-edges and thus $U$ is an orbit of $\Aut(\Gamma)$, contradicting the assumption that $\Gamma$ is vertex-transitive. Hence, $\Gamma_{0,R}$ is disconnected.
\end{proof}

\color{black}

\begin{lemma}
\label{lem:4.6}
If $3r-s \equiv 0$ $(\mod 2k)$, $k \geq 9$, then $T_1(k,r,s)$ is not vertex-transitive. 
\end{lemma}

\begin{proof}
Let $\Gamma=T_1(k,r,s)$ where $3r-s \equiv 0$ $(\mod 2k)$. A computer assisted counting shows that $S$-edges and $K$-edge have exactly $4$ distinct $8$-cycles passing through them, while $0$- and $R$-edges have $6$ such cycles (see Table \ref{table:2}). It follows that any automorphism sending a vertex $u \in U$ to a vertex $v \in V$ must necessarily map the $K$-edge incident to $u$ into the $S$-edge incident to $v$. Thus $E_S \cup E_K$ is an orbit of edges under the action of $\Aut(\Gamma)$.

Since we are assuming that $3r-s \equiv 0$ we get that any number dividing $k$ and $r$ must also divide $s$, but since $\gcd(k,r,s)=1$, we have that $\gcd(k,r)=1$ and then $\gcd(2k,r) \in \{1,2\}$. In light of Lemma \ref{lem_blocks}, $\Gamma_{0,R}$ is disconnected, implying $\gcd(2k,r) \neq 1$ and therefore $\gcd(2k,r) = 2$. From the congruence $3r-s \equiv 0$ we also get that $r$ and $s$ must have the same parity thus making $s$ even and $k$ odd. Recall that each connected component of $\Gamma_{0,R}$ consists of all the vertices whose indices are congruent modulo $\gcd(2k,r)$, so $\Gamma_{0,R}$ has two connected components: one containing all the vertices with even index, and the other containing those with odd index. Since $s$ is even and $k$ is odd, $S$-edges in $\Gamma$ join vertices with same parity, while $K$-edges join vertices with distinct parity. We know that each of the two connected components of $\Gamma_{0,R}$ is a block of imprimitivity of $\Gamma$, an so any automorphism of $\Gamma$ must either preserve the parity of all indices, or of none at all. It follows that no automorphism can send $S$-edges into $K$-edges, and therefore no automorphism can send a vertex in $V$ to a vertex in $U$. We conclude $\Gamma$ cannot be vertex-transitive. 
\end{proof}

\begin{lemma}
\label{lem:4.7}
If $3s-r \equiv 0$  $(\mod 2k)$, $k \geq 9$, then $T_1(k,r,s)$ is not vertex-transitive.
\end{lemma}

\begin{proof}
Set $\Gamma=T_1(k,r,s)$ and $\Gamma'=T_1(k,r',s')$, where $r'=s$ and $s'=r$. Notice that the triplet $(k,r',s')$ satisfies equation \ref{eq3} and so $\Gamma'$ cannot be vertex-transitive. By observation \ref{eq:rs}, $\Gamma \cong \Gamma'$. Therefore $\Gamma$ is not vertex-transitive.
\end{proof}

\begin{lemma}
\label{lem:4.8}
If $4r-4s \equiv 0$ $(\mod 2k)$, $k \geq 9$, then $T_1(k,r,s)$ is not vertex-transitive.
\end{lemma}

\begin{proof}
Let $\Gamma=T_1(k,r,s)$, with $4r-4s \equiv 0$ $(\mod 2k)$ and $k \geq 9$. From the congruence $4r-4s \equiv 0$, we see that either $2r-2s \equiv 0$ or $2r - 2s + k \equiv 0$, but the former contradicts Lemma \ref{lem:girth5}. Hence $2r - 2s + k \equiv 0$ and the closed walk $((uv)_0,(vw)_r,(wv)_{-s},(vw)_r,(wv)_{-s},(vu)_0,(uu)_k)$ in $\Delta_1$ lifts to a $7$-cycle in $\Gamma$. We will show that $\Gamma$ cannot be $7$-vertex-regular. Following a similar procedure as the one used to obtain Table \ref{table:1}, we have counted all closed walks of length $7$ in $\Delta_1$ starting at $u$ and those starting at $v$, and we have computed their net voltage. The result is displayed in Table \ref{table:3}.

\begin{table}
\begin{tabular}{|c| c| c| c|}
\hline
 Label & Net voltage & Starting at $u$ & Starting at $v$   \\
 \hline
 $\textrm{I}$ & $\pm (k+2r-2s)$ & $8$ & $10$ \\
 $\textrm{II}$ & $\pm (k+r+s)$ & $12$ & $8$  \\
 $\textrm{III}$ & $\pm (k+2r)$ & $6$ & $4$  \\
 $\textrm{IV}$ & $\pm (k+2s)$ & $6$ & $4$  \\
 $\textrm{V}$ & $\pm (3r-2s)$ & $2$ & $6$ \\
 $\textrm{VI}$ & $\pm (3s-2r)$ & $2$ & $6$  \\
 \hline
\end{tabular}
\caption{Net voltages of closed walks of length $7$ in $\Delta_1$.}
\label{table:3}
\end{table}

It is plain to see that if, in addition to $\textrm{I}$, any one of the voltages in rows $\textrm{III}$--$\textrm{VI}$ is congruent to $0$, then either $r=0$, $r=k$, $s=0$ or $s=k$, contradicting Lemmas \ref{lem_rsnok} or \ref{lem:rsno0}. We see that only the voltages in row $\textrm{I}$ or $\textrm{II}$ can hold, and thus $\Gamma$ is not $7$-vertex-regular and hence cannot be vertex-transitive. 
\end{proof}

In what follow we deal with the case where (\ref{eq1}) holds, that is, when $3s-2r+k \equiv 0$ (\mod $2k$) and $k\geq9$. From Table \ref{table:2} we see that $S$-edges and $K$-edge have exactly $6$ distinct $8$-cycles passing through them, while $0$- and $R$-edges have only $5$. In particular, no automorphism of $\Gamma$ maps an $S$-edge to a $0$, implying that $E_S \cup E_K$ is an edge-orbit of $\Gamma$ (see the discussion above Lemma~\ref{lem_blocks}).

\begin{lemma}
\label{lem:ks}
Suppose that  $3s-2r+k \equiv 0$ (\mod $2k$), $k\geq9$.
If $e$ is an edge in $E_K \cup E_S$, then the endpoints of $e$ belong to different connected components of $\Gamma_{0,R}$.
\end{lemma}

\begin{proof}
Suppose a $K$-edge $e$ has both its endpoints in the same connected component, $C$, of $\Gamma_{0,R}$. Then, because $\Gamma$ is vertex-transitive and $C$ is a block of imprimitivity, both of the endpoints of any $K$-edge must belong to the same connected component of $\Gamma_{0,R}$. This means that the subgraph induced by $0$-edges, $R$-edges and $K$-edges is disconnected. Recall that $K$-edges and $S$-edges belong to a single edge-orbit in $\Gamma$ and thus, for every $S$-edge, $e'$, there exists $\phi \in \Aut(\Gamma)$ such $\phi(e)=e'$. Since $\phi$ also acts as an automorphism on $\Gamma_{0,R}$, it follows that the endpoints of $\phi(e)$ are contained in the component $\phi(C)$, thus making $\Gamma$ a disconnected graph, which is contradiction. 
\end{proof}

\begin{lemma}
\label{lem:evencc}
Suppose that  $3s-2r+k \equiv 0\> (\mod 2k)$, $k\geq9$. Then
the subgraph $\Gamma_{0,R}$ has an even number of connected components. 
\end{lemma}

\begin{proof}
Let $e$ be a $K$-edge whose endpoints are in two different connected components $C_1$ and $C_2$. Let $\rho^r \in \Aut(\Gamma)$ be an "$r$-fold" rotation. That is, $\rho^r$ adds $r$ to the index of every vertex. It is plain to see that $\rho^r$ fixes the connected components of $\Gamma_{0,R}$ set-wise and that it send $K$-edges into $K$-edges. Moreover, since it is transitive on $U \cap C_1$, we have that all $K$-edges having an endpoint in $C_1$ will have the other endpoint in $C_2$. Hence $K$-edges 'pair up' connected components of $\Gamma_{0,R}$. The result follows
\end{proof}

\begin{figure}[h!]
\centering
\includegraphics[width=0.8\textwidth]{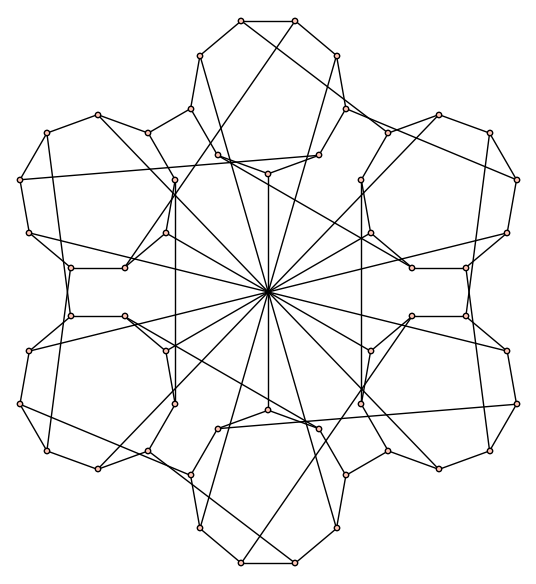}
\caption{$T_1(9,6,1)$}
\label{fig:t961}
 \end{figure}

\begin{proposition}
\label{prop:necessary}
If $T_1(k,r,s)$ is connected and vertex-transitive, then the following statements (or the analogous three statements obtained by interchagning $r$ and $s$) hold:
\begin{enumerate}
\item $3s-2r+k \equiv 0$ $( \mod 2k )$
\item $k$ and $s$ are odd and $\gcd(k,s)=1$
\item $r$ is even and $\gcd(k,r)=\{1,3\}$
\end{enumerate}
\end{proposition}

\begin{proof}

First, note that item (1) follows from the combination of Lemma \ref{lem:eqtype1} and Lemmas \ref{lem:4.6}, \ref{lem:4.7} and \ref{lem:4.8}.

We will continue by proving (3). Set $d=\gcd(k,r)$, so that $r=dr'$ and $k=dk'$ for some integers $k'$ and $r'$. Since $3s-2r+k \equiv 0\> (\mod 2k)$, we see that, for some odd integer $a$:
\begin{eqnarray*}
3s &=& d(ak'+2r')
\end{eqnarray*}

By connectedness, $\gcd(k,r,s)=1$, and hence $\gcd(d,s)=1$ implying that $d$ divides $3$. Hence $d = \gcd(k,r)\in \{1,3\}$. Moreover, since $\Gamma_{0,R}$ has $\gcd(2k,r)$ connected components and this number must be even by Lemma \ref{lem:evencc}, we see that $\gcd(2k,r)\in \{2,6\}$, and therefore $r$ is even.

 Now, from the congruence $3s-2r+k \equiv 0\> (\mod 2k)$ we see that $s$ and $k$ are either both even, or both odd. Since $\gcd(2k,r,s)=1$ and $r$ is even, $s$ must necessarily be odd and then $k$ is also odd. Set $d':=\gcd(k,s)$ so that $s=d's'$ and $k=d'k'$ for some integers $s'$ and $k'$. Again, from equation (\ref{eq1}) we obtain
 $$2r=d'(ak'+ 3s')$$
 for some odd integer $a$. Thus $d' \in \{1,2 \}$, but since $s$ is odd, we see that $\gcd(k,s)=1$. This proves item (2) and completes the proof.
\end{proof}

\begin{remark}
\label{rem:type1}
Conditions (1)--(3) (or their analogues) in Proposition \ref{prop:necessary} are also sufficient to prove vertex-transitivity of a connected $T_1(k,r,s)$. This will be shown after the proof on Lemma \ref{lem:isot1}. Moreover, observe that if $k$ and $s$ satisfy condition (2) of Lemma \ref{prop:necessary}, then $\gcd(2k,1)=1$ and $T_1(k,r,s)$ is isomorphic to $T_1(k,rs^{-1},1)$ where $s^{-1}$ is the multiplicative inverse of $s$ in $\ZZ_{2k}$. Hence, in what follows we may assume that $s=1$.
\end{remark}

\begin{lemma}
\label{lem:detr}
Let $k$ be an odd number and set $s=1$. There exists a unique element $r \in \ZZ_{2k}$ such that conditions (1) and (3) of Proposition \ref{prop:necessary} are satisfied. This unique $r$ equals $r^*$ from Definition \ref{def:X}.
\end{lemma} 

\begin{proof}
Suppose that such an $r$ exists. From (1), we get
$r = (3-ak)/2$
for some odd integer $a$. Since $r$ must be a positive integer smaller than $2k$, it follows that $a \in \{-1,-3\}$. Hence $(3+k)/2$ and $(3+3k)/2=(3+k)/2+k$ are the only possible candidates for $r$. Since $k$ is odd, the integers $(3+k)/2$ and $(3+k)/2+k$ have different parity. Set $r$ to be whichever one of these two expressions is even. Note that any odd number that divides both $r$ and $k$ must divide $3$ and that this is true whether $r$ equals $(3+k)/2$ or $(3+k)/2+k$.  Thus conditions (3) is satisfied and $r$ equals $r^*$ from Definition \ref{def:X}.
\end{proof}

We have just proved that once an odd $k \geq 9$ is prescribed, then there is at most one graph $T_1(k,r,1)$ that is connected and vertex-transitive. The following lemma generalizes this to an arbitrary value of $s$.

\begin{lemma}
\label{lem:isot1}
Let $k \geq 9$ be an odd integer and let $T_1(k,r,s)$ be vertex-transitive. Then $T_1(k,r,s) \cong \textup{X}(k)$.
\end{lemma} 

\begin{proof}
Recall that $\textup{X}(k) \cong T_1(k,r^*,1)$, where $r^*$ is as in Definition \ref{def:X}. Since $T_1(k,r,s)$ is connected and vertex-transitive, conditions (1)--(3) of Proposition \ref{prop:necessary} hold. In particular, $s$ is relatively prime to $2k$. Then  $T_1(k,r,s) \cong T_1(k,s^{-1}r,1)$, where $s^{-1}$ is the multiplicative inverse of $s$ in $\ZZ_{2k}$. By Lemma~\ref{lem:detr}, $T_1(k,s^{-1}r,1) = T_1(k,r^*,1)$.
\end{proof}

Notice that in the previous proof vertex-transitivity is only used to ensure that the graph $T_1(k,r,s)$ satisfies conditions (1)--(3) of Proposition \ref{prop:necessary}, and from there $T_1(k,r,s)$ is shown to be isomorphic to $T_1(k,r^*,1)$. Since we know from Lemma~\ref{lem: xvt} that $T_1(k,r^*,1)$ is vertex-transitive, we have that conditions (1)--(3) of Proposition~\ref{prop:necessary} are not only necessary, but also sufficient for vertex-transitivity, and thus Proposition~\ref{prop:necessary} may be regarded a characterization theorem.

Now that we have a chracterization for vertex-transitivity, we would like to know when a cubic vertex-transitive tricirculant of Type $1$ is also a bicirculant. Note that the only cubic vertex-transitive circulants are prism and M\"obius ladders, which have girth $4$. It follows from Lemma \ref{lem:girth5} that no cubic vertex-transitive tricirculant of Type $1$ is a circulant.

\begin{lemma}
\label{lem:t1bic}
Let $\Gamma = T_1(k,r,s)$ be connected and vertex-transitive. If $3 \nmid k$, then $\Gamma$ is a bicirculant.
\end{lemma}

\begin{proof}
Let $\varphi: V(\Gamma) \to V(\Gamma)$ be the mapping given by:
\begin{center}
\begin{tabular}{llll}
$u_i \mapsto v_i$, & $v_i \mapsto w_{i+r}$, & $w_i \mapsto u_i$ & if $i$ is even;  \cr
$u_i \mapsto w_{i+k+s}$, & $v_i \mapsto u_{i+k+s}$, & $w_i \mapsto v_{i+k+s-r}$ & if $i$ is odd. 
\end{tabular}
\end{center} 

It can be readily seen that $\varphi$ is indeed a graph automorphism. Moreover, since $r$ is even and both $k$ and $s$ are odd, $\varphi$ preserves the parity of the index of all vertices. In fact, it is plain to see that $\varphi$ is the product of the following two disjoint permutation cycles $\phi_1$ and $\phi_2$ of length $3k$:
\begin{eqnarray*}
\phi_0 &=&(u_0,v_0,w_r,u_r,v_r,w_{2r},u_{2r},v_{2r},...,u_{(k-1)r},v_{(k-1)r},w_0); \cr
\phi_1 &=&(u_k,w_{s},v_{k+2s-r},u_{k+r},w_{r+s},v_{k+2s},u_{k+2r}...v_{k+2s+(k-2)r}).   
\end{eqnarray*}
Hence, $\Gamma$ is a bicirculant.
\end{proof}

\begin{corollary}
\label{cor:GP}
If $T_1(k,r,s)$ is vertex-transitive with $3 \nmid k$, then $T_1(k,r,s) \cong \hbox{{\rm GP}}(3k,k+(-1)^\alpha)$ where $\alpha \in \{1,2\}$ and $\alpha \equiv k$ (\mod$3$).
\end{corollary}

\begin{proof}
Let $T_1(k,r,s)$ be vertex-transitive and suppose $k \equiv 1$ (\mod $3$). Let $\varphi=\phi_0\phi_1$ the automorphism described in the proof of Lemma \ref{lem:t1bic}. 
 We will show that $u_k$ is adjacent to $\varphi^{(k+1)}(u_k)$. Notice that $\varphi^{(k+1)}(u_k)=v_{k+2s-r+\frac{k-1}{3}r}$, but $k+2s-r \equiv r-s$ so we may write $\varphi^{(k+1)}(u_k)=v_{r-s+\frac{k-1}{3}r}$. Recall that $\gcd(2k,3)=1$ and $r$ is even. From the congruence $k+2r-3s \equiv 0$ we see that:
\begin{eqnarray*}
2r-3s &\equiv& 3k \cr
3r-3s+r(k-1) &\equiv& 3k \cr
r-s+r(k-1)/3 &\equiv& k 
\end{eqnarray*}
So that $\varphi^{(k+1)}(u_k)=v_k$, which is adjacent to $u_k$. The case when $k \equiv 2$ follows from a similar argument. 
\end{proof}

We have now completed the proof of Theorem \ref{the:t1}. The characterization of vertex-transitivity for cubic tricirculants of Type $1$ with at least $54$ vertices follows from Lemmas \ref{lem: xvt} and \ref{lem:isot1}. That a vertex-transitive tricirculant of Type $1$ is a Generalized Petersen graph when $3 \nmid k$ follows from Lemma \ref{lem:t1bic} and Corollary \ref{cor:GP}.

\section{Type 2} \label{Type2}

\begin{figure}[h!]
\centering
\includegraphics[width=0.3\textwidth]{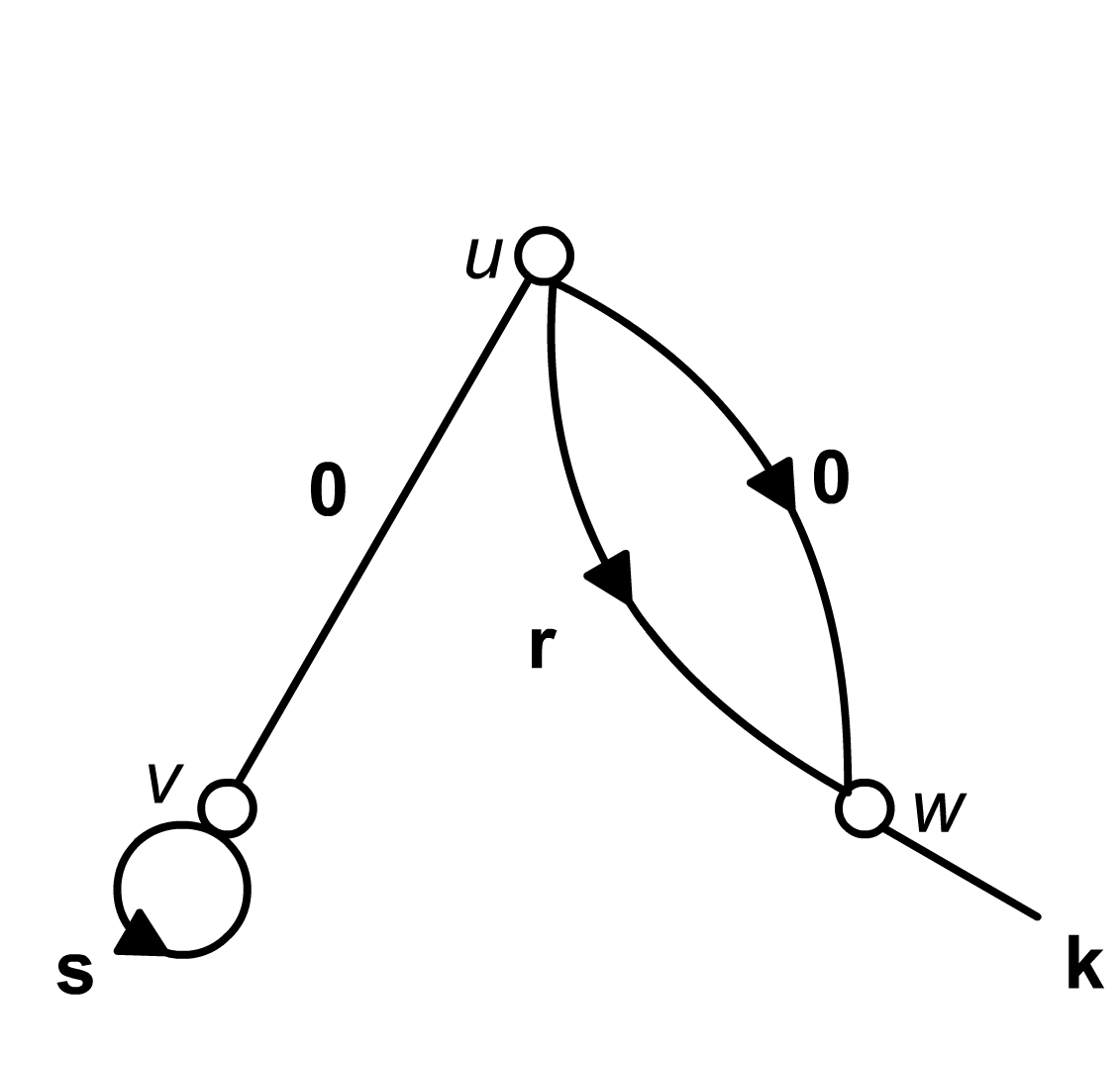}
\caption{Type 2}
\label{fig:Delta2}
 \end{figure}

Let $k \geq 9$ be an integer, and let $r,s \in \ZZ_{2k}$. Define $T_2(k,r,s)$ as the derived graph of $\Delta_2$ with the normalized voltage assignment for $\ZZ_{2k}$ shown in Figure \ref{fig:Delta2}.\\

Then, with the same notation as in Section \ref{type1}, $U=\{u_0,u_1,...,u_{2k-1}\}$, $V=\{v_0,v_1,...,v_{2k-1}\}$ and $W=\{w_0,w_1,...,w_{2k-1}\}$ are the respective fibres of vertices $u$, $v$ and $w$ in $\Delta_2$. The set of edges of $T_2(k,r,s)$ can be expressed as the union $E_K \cup E_R \cup E_S \cup E_0$ where:\\
\begin{eqnarray}
E_K &=& \{w_iwu_{i+k}: i \in \ZZ_{2k}\} \cr
E_0 &=& \{u_0v_0:i \in \ZZ_{2k}\} \cup \{u_0w_0:i \in \ZZ_{2k}\} \cr
E_R &=& \{u_iw_{i+r}: i \in \ZZ_{2k}\} \cr 
E_S &=& \{v_iv_{i+s}: i \in \ZZ_{2k}\} \cr
\end{eqnarray}
Similarly as with Type $1$ cubic tricirculants, we see that every cubic tricirculant of Type $2$ is isomorphic to $T_2(k,r,s)$ for an appropriate choice of $k$ and $r,s \in \ZZ_{2k}$.

\begin{definition}
\label{def:Y}
For an odd positive integer $k$, let $\textup{Y}(k)= T_2(k,2,1)$.
\end{definition}

\begin{theorem}
\label{theo:t2}
Let $\Gamma$ be a connected cubic tricirculant of Type 2 with at least $54$ vertices.
Then $\Gamma$ is vertex-transitive if and only if it is isomorphic to $\textup{Y}(k)$ for some odd integer $k$, $k\ge 9$.
In this case, $\Gamma$ is a bicirculant if $3 \nmid k$. Furthermore, $\Gamma$ is a map on the torus of type $\{6,3\}_{\alpha,3}$ where $\alpha = \frac{1}{2}(k-3)$.
\end{theorem}

\begin{remark}
\label{rem:-ss}
Note that $T_2(k,r,s)$ and $T_2(k,r,-s)$ are in fact the exact same graph. Then, we can safely assume $s<k$.
\end{remark}

\begin{lemma}
\label{lem:yk}
For an odd integer $k$, $k \geq 9$, $\textup{Y}(k)$ is vertex-transitive.
\end{lemma} 

\begin{proof}
Since $\textup{Y}=T_2(k,r,s)$, it suffices to provide an automorphism of $T_2(k,2,1)$ that mixes the sets $U$, $V$ and $W$. Let $\varphi$ be the mapping given by:

\begin{center}
\begin{tabular}{llll}

$u_i \mapsto v_{i+1}$, & $v_i \mapsto u_{i+1}$, & $w_i \mapsto v_i$ & if $i$ is even; \cr
$u_i \mapsto w_{i+2+k}$, & $v_i \mapsto w_{i+2}$, & $w_i \mapsto u_{i+k}$ & if $i$ is odd.

\end{tabular}
\end{center}

Now consider a vertex $u_i\in U$ with $i$ even. Observe that $\varphi$ maps $u_i$ to $v_{i+1}$ and that the neighbourhood $N(u_i):=\{v_i,w_i,w_{i+2}\}$ is mapped to $\{u_{i+1},v_i,v_{i+2}\}$, which is precisely the neighbourhood of $v_{1+i}$. That is, $\varphi$ maps the neighbourhood of any vertex $u_i$ into the neighbourhood of its image, when $i$ even. The reader can verify the remaining cases and see that $\varphi$ is indeed a graph automorphism.
\end{proof}

For the rest of this section, let $k \geq 9$ be an integer, let $r,s \in \ZZ_{2k}$ and suppose $\Gamma=T_2(k,r,s)$ is vertex-transitive. In what follows we will show that $\Gamma$ is isomorphic to $\textup{Y}(k)$. As with cubic tricirculants of Type 1, the strategy will be to count reduced closed walks in the quotient $\Delta_2$. Observe that the walk $((wu)_0,(uw)_r,(ww)_k,(wu)_{-r},(uw)_0,(ww)_k)$ starting at $w$ is a closed walk of length $6$ having net voltage $0$, regardless of what $r$ and $s$ evaluate to in $\ZZ_{2k}$. As was done with graphs of Type 1, we computed all closed walks in $\Delta_2$ along with their net voltages. Again, it would be convenient if closed walks of length $6$ with net voltage $0$ always lifted into $6$-cycles. For this it suffices to show that the girth of $\Gamma$ is at least $4$.

\begin{lemma}
\label{lem:triangles}
$\Gamma$ has no cycles of length $3$.
\end{lemma}

\begin{proof}
Suppose to the contrary, that $\Gamma$ contains a triangle $T$. By observing Figure \ref{fig:Delta2} we know that all three vertices of $T$ belong to $V$ or they belong to $U \cup W$. Since $\Gamma$ is vertex-transitive, there is at least one triangle through every vertex in $W$. We can thus assume without loss of generality that all three vertices of $T$ are in $U \cup W$. This implies that $r=k$ and further, that there are $2$ distinct triangles through every vertex in $W$: the lifts of $((wu)_0,(uw)_r,(ww)_k)$ and $((wu)_{-r},(uw)_0,(ww)_k)$. However, since the subgraph induced by $V$ is $2$-valent, there is at most one triangle through every vertex in $V$, contradicting the vertex-transitivity of $\Gamma$.
\end{proof}

Note that since $\Gamma$ is simple, $s\neq k$ and if $r=k$, $\Gamma$ would contain a triangle. Thus

\begin{corollary}
Neither $r$ nor $s$ equals $k$.
\end{corollary}

Now, let $\mW_6$ be set of all walks of length $6$ in $\Delta_2$ and let $N$ be the set of net voltages of elements of $\mW_6$, expressed in terms of $k$, $r$ and $s$. It follows from Lemma \ref{lem:triangles} that walks in $\mW_6$ with net voltage $0$ will lift into $6$-cycles, and not just closed walks. Table \ref{table:4} shows, for each element $n$ of $N$, how many closed walks with  net voltage $n$ start at each vertex of $\Delta_2$.

\begin{center}
\begin{table}
\begin{tabular}{|c| c| c| c| c|}
\hline
  Label & Net voltage & Starting at $u$ & Starting at $v$ & Starting at $w$  \\
 \hline
 \textrm{I} & $0$ & $2$ & $0$ & $4$ \\
 \textrm{II} & $\pm(k-s)$ & $8$ & $8$ & $8$  \\
 \textrm{III} & $\pm(k-r-s)$ & $4$ & $4$ & $4$  \\
 \textrm{IV} & $\pm(k+r-s)$ & $4$ & $4$ & $4$  \\
 \textrm{V} & $\pm(3r)$ & $2$ & $0$ & $2$  \\
 \textrm{VI} & $\pm(2r)$ & $2$ & $0$ & $4$  \\
 \textrm{VII} & $\pm(6s)$ & $0$ & $2$ & $0$  \\
 \textrm{VII} & $\pm(r-2s)$ & $4$ & $6$ & $2$ \\
 \textrm{IX} & $\pm(r+2s)$ & $4$ & $6$ & $2$  \\
  \hline 
\end{tabular}
\caption{Net voltages of closed walks of length $6$ in $\Delta_2$.}
\end{table}
\label{table:4}
\end{center}

From row I of Table \ref{table:4} we see that there are at least $4$ walks $\omega \in \mW_6$ starting at $w$. Therefore, at least one expression from rows II--IX must be congruent to $0$. However, a careful inspection of Table \ref{table:4} shows that if neither of the expressions in rows VIII and IX are congruent to $0$, then there are at least $2$ more walks in $\mW_6$ with net voltage $0$ starting at $w$ than there are those starting at $v$. This means that if $\Gamma$ is vertex-transitive, then one of the two following equations modulo $2k$ hold:
\begin{eqnarray}
\label{eqA} r \equiv 2s\\
\label{eqB} r \equiv -2s
\end{eqnarray}
Suppose \ref{eqA} holds. By Remark \ref{rem:-ss} we can in fact write $r = 2s$. Then the mapping $\phi: T_2(k,2s,s) \to T_2(k,-2s,-s)$ given by $x_i \mapsto x_{-i}$, for all $x \in U \cup V \cup W$, is a graph isomorphism. But, again by Remark \ref{rem:-ss}, $T_2(k,-2s,-s) = T_2(k,-2s,s)$ so that any graph $T_2(k,r,s)$ with $r \equiv 2s$  is isomorphic to a graph $T_2(k,r',s)$ satisfying $r' \equiv -2s$. We can therefore limit our analysis to the case when \ref{eqA} holds.

Let $\Gamma = T_2(k,r,s)$ be connected and vertex-transitive, and suppose $r \equiv 2s$. Note that any number dividing both $k$ and $s$ must also divide $r$, and since by connectedness of $\Gamma$ $\gcd(k,r,s)=1$ we see that $\gcd(k,s)=1$ and $\gcd(2k,s) \in \{1,2\}$. We will show that $\gcd(2k,s)$ cannot be $2$.

\begin{lemma}
If $\Gamma$ is vertex-transitive, $\gcd(2k,s)=1$ and $\gcd(2k,r)=2$.
\end{lemma}

\begin{proof}
In order to get a contradiction, suppose $\gcd(2k,s)=2$. Then $s$ is even and $k$ is odd, as $\gcd(k,s)=1$. Furthermore, the $2$-valent subgraph of $\Gamma$ induced by the set $V$ is the union of two $k$-cycles. We will show that these are in fact the only two $k$-cycles in $\Gamma$ which will imply that $\Gamma$ cannot be vertex-transitive. Suppose there is a $k$-cycle $C'$ that visits a vertex not in $V$. Define $\Gamma[U \cup W]$ as the subgraph of $\Gamma$ induced by the set $U \cup W$. Since $r$ is even and $k$ is odd, it is not difficult to see that $\Gamma[U \cup W]$ is bipartite, with sets $\{w_i:\text{$i$ is even}\} \cup \{u_i:\text{$i$ is odd}\}$ and $\{w_i:\text{$i$ is odd}\} \cup \{u_i:\text{$i$ is even}\}$. Therefore no $k$-cycle of $\Gamma$ is contained in  $\Gamma[U \cup W]$. This means that $C'$ visits at least one vertex in each of sets $U$, $V$ and $W$. 

Now, $C'$ projects onto a closed walk $C$ in $\Delta_2$ that has net voltage $0$. Define $x_S$ as the number of times $C$ traces the dart $(vv)_s$ minus the number of times it traces $(vv)_{-s}$, so that $x_Ss$ is the total voltage contributed to $C$ by darts in $\{(vv)_{\pm s}\}$. Define $x_R$ similarly, and define $x_K$ as the number of times $C$ traces the semi-edge $(ww)_k$. We obtain the following congruence modulo $2k$:
$$ x_Rr + x_Ss + x_Kk \equiv 0.$$
Recall that both $r$ and $s$ are even, making both $x_Rr$ and $x_Ss$ even. It follows that $x_Kk$ is even, but since $k$ is odd, $x_K$ must be even, and hence $x_Kk \equiv 0$ $(\mod 2k)$. We thus have 
$$ x_Rr + x_Ss \equiv 0$$
and since $r \equiv 2s$ $(\mod2k)$, 
\begin{eqnarray*}
x_R2s + x_Ss \equiv 0, \cr
s(2x_R + x_S) \equiv 0.
\end{eqnarray*}
From this, and because $\gcd(k,s)=1$, we see that $2x_R + x_S = 0$ or $2x_R + x_S \geq k$. To see that $2x_R + x_S$ cannot equal $0$, define $A$ as the set of darts in $\Delta_2$ with voltage $0$ or $r$ and observe that $C$ must visit $A$ an even number of times. Since $C$ has odd length, $C$ visits $\{(vv)_{\pm s}\} \cup \{(ww)_k\}$ an odd amount of times. Then, since $x_K$ is even, $x_S$ must be odd. But $2x_R$ is even so $2x_R + x_S \neq 0$. We thus have
\begin{equation*}
\label{geqk}
2x_R + x_S \geq k. 
\end{equation*}
Now, notice that if $C$ traces the dart $(uw)_{r}$, then it must immediately trace a dart in $\{(wu)_0,(ww)_k\}$. This means that $C$ traces a dart in the subgraph of $\Delta_2$ induced by $U \cup W$ at least $2x_R$ times. Then $2x_R + x_S < k$, since $C$ has length $k$ and it must visit $(vu)_0$ at least once. We thus have $k > 2x_R + x_S  \geq k$, a contradiction. The result follows.
\end{proof}

\begin{lemma}
If $T_2(k,r,s)$ is vertex-transitive, $k$ is odd.
\end{lemma}

\begin{proof}
Suppose to the contrary that $k$ is even. Since $s$ is odd, we have $sk \equiv k$ $(\mod 2k)$, and thus $s(2\cdot\frac{1}{2})k \equiv k$. But $k$ is even, so we can rewrite that expression as $2s\cdot\frac{k}{2} \equiv k$. Recall that $r \equiv 2s$ and then $r\frac{k}{2} \equiv k$. 

Now, consider the walk of length $k+1$ in $\Delta_2$ that starts in $w$, traces the semi-edge $(ww)_k$ once and then traces the $2$-path $((wu)_r,(uw)_0)$ exactly $\frac{k}{2}$ times. Note that this walk is closed and has net voltage $\frac{k}{2}r + k \equiv 0$. Furthermore, it is easy to see that it lifts into a $(k+1)$-cycle in $\Gamma$ and that it doesn't visit any vertex in $V$.

Since $\Gamma$ is vertex-transitive, there must be a $(k+1)$-cycle $C'$ through each vertex in $V$. Note that the subgraph of $\Gamma$ induced by $V$ is a single $2k$-cycle, and thus $C'$ must visit all three set $U$, $V$ and $W$. It follows that $C'$ projects into a walk $C$ of length $k+1$ having net voltage $0$ and that it visits all three vertices of $\Delta_2$. Furthermore, since $k+1$ is odd, the number of times $C$ visits $\{(vv)_{\pm s}\}$ must be of different parity than the number of times it traces $(ww)_k$. Define $x_S$, $x_R$ and $x_K$ as in the proof of Lemma 19. Hence
$$x_Rr + x_Ss + x_Kk \equiv 0.$$
Now, $r$ and $k$ are even, and so $x_Rr$ and $x_Kk$ are also even. It follows that $x_Ss$ is even, but since $s$ is odd, $s_S$ is necessarily even. This means  $C$ visits  $\{(vv)_{\pm s}\}$ an even amount of times. It follows that $C$ traces $(ww)_k$ and odd amount of times, that is, $x_K$ is odd. We thus have
\begin{eqnarray*}
x_Rr + x_Ss + k \equiv 0, \cr
x_Rr + x_Ss \equiv k,
\end{eqnarray*}
and since $r \equiv 2s$,
$$s(2x_R + x_S) \equiv k.$$
This implies $2x_R + x_S \geq k$, as $s$ and $k$ are relatively prime. However, $2x_R + x_R \leq k-1$ since $C$ has length $k+1$ and $C$ visits the set $\{(vu)_0,(uv)_0\}$ at least twice. This contradiction arises from the assumption that $k$ is even. We conclude that $k$ is odd. 
\end{proof}
 
\begin{lemma}
\label{lem:t2cong}
For an odd $k>9$, if $T_2(k,r,s)$ is vertex-transitive then $T_2(k,r,s) \cong T_2(k,2,1)$.
\end{lemma}

\begin{proof}
Let $T_2(k,r,s)$ be vertex-transitive. Then, $k$ is an odd integer, $\gcd(2k,s)=1$ and $r \equiv 2s$. Denote by $s^{-1}$ the multiplicative inverse of $s$ in $\ZZ_{2k}$. Now, define $\varphi$ as the mapping between $T_2(k,2,1)$ and $T_2(k,r,s)$ given by $x_i \mapsto x_{s^{-1}i}$, for $x \in U \cup V \cup W$. It is plain to see that $\varphi$ is the desired isomorphism.
\end{proof}

This proves the first claim of Theorem \ref{theo:t2}. We now proceed to show sufficient conditions under which $\Gamma$ is a bicirculant and we show that $\Gamma$ can be seen as a map on the torus.

\begin{lemma}
\label{lem:t2bic}
If $k \geq 9$ is an odd integer such that $3 \nmid k$, then $T_2(k,2,1)$ is a bicirculant.
\end{lemma}

\begin{proof}
Consider the mapping $\varphi$ defined by:

\begin{center}
\begin{tabular}{llll}
$u_i \mapsto v_{i+1}$, & $v_i \mapsto u_{i+1}$, & $w_i \mapsto v_i$, & if $i$ is even; \cr
$u_i \mapsto w_{i+2+k}$, & $v_i \mapsto w_{i+2}$, & $w_i \mapsto u_{i+k}$, & if $i$ is odd.
\end{tabular}
\end{center}

We have that for a vertex $u_i \in U$, $\varphi^l(u_i) \in U$ if and only if $3 | l$. This is, if we start at $U$, every third iteration of $\varphi$ lands us in $U$ again. Moreover, we have $\varphi^3(u_i)=u_{i+3+k}$. In general $\varphi^{3l}(u_i)=u_{i+3l+lk}$. Hence $\varphi^{3l}(u_i)=u_i$ if and only if $3l + lk \equiv 0$. If $3 \nmid k$, then $l = k$ is the smallest value for $l$ that satisfies $3l + lk \equiv 0$. It follows that the orbit of $u_i$ has size $3k$. It is plain to see that the vertices not in this orbit form an orbit of size $3k$ on their own, namely, the orbit of $u_{i+1}$. It is worthwhile to mention that if $3$ does divide $k$, then $\varphi$ has $6$ orbits of size $k$.         
\end{proof}

\begin{proposition}
\label{prop:torus}
Let $k\geq 3$ be an odd integer, then $T_2(k,2,1)$ admits an embedding on the torus with hexagonal faces, yielding a map of type $\{6,3\}_{\alpha,3}$, where $\alpha = \frac{1}{2}(k-3)$ (see \cite{torus} for details about the notation).
\end{proposition}

\begin{figure}[h!]
\centering
\includegraphics[width=0.8\textwidth]{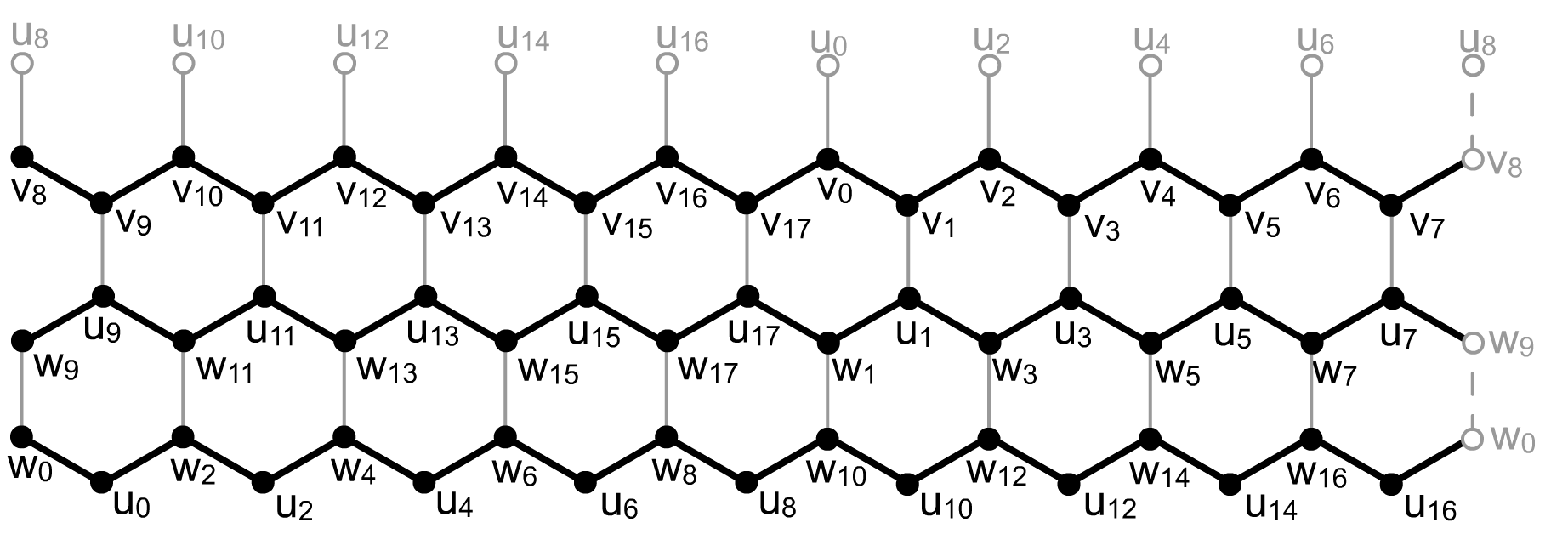}
\caption{$T_2(9,2,1)$. The subgraph induced by bold edges has three connected components that correspond, from bottom to top, to the cycles $C_1$, $C_2$ and $C_3$ described in the proof of Proposition \ref{prop:torus}.}
\label{fig:torus1}
\end{figure}

\begin{proof}
First, notice that vertices in $U \cup W$ with  even index form a cycle of length $2k$, $C_1=(w_0,u_0,w_2,u_2,...,w_{2k-2},u_{2k-2})$. Likewise, vertices  in $U \cup W$ with odd index form a $2k$-cycle, $C_2$. The vertices in $V$ form a third cycle of length $2k$, $C_3$.

Now, $w_iw_{i+k} \in E$ for all $i \in \ZZ_{2k}$, and $i+k$ has different parity than $i$. This means each vertex of $C_1$ of the form $w_i$ is adjacent, through a $K$-edge, to a vertex of $C_2$. Further, each vertex of the form $u_i$ in $C_1$ is adjacent to a vertex in $C_3$, namely $v_i$. Note that the vertices of $C_2$ that are not adjacent to a vertex of $C_1$ are precisely those of the form $u_i$ (with $i$ odd), and that $u_iv_i \in E$ for all odd $i$. This is, every other vertex in $C_2$ has a neighbour in $C_3$.

Therefore, we can think of $T_2(k,2,1)$ as three stacked $2k$-cycles $C_1$, $C_2$ and $C_3$ where every other vertex of $C_i$ has a neighbour in $C_{i-1}$ while each of the remaining vertices has a neighbour in $C_{i+1}$, where $i \pm 1$ is computed modulo $3$ (see Figure \ref{fig:torus1} for a detailed example). From here it is clear that $T_2(k,2,1)$ can be embedded on a torus, tessellating it with $3k$ hexagons. It can be readily verified that this wields the map $\{6,3\}_{\alpha,3}$, where $\alpha = \frac{1}{2}(k-3)$, following the notation in \cite{torus}.
\end{proof}

Observe that $T_2(9,2,1)$ corresponds to the map $\{6,3\}_{3,3}$, which is a regular map (see Chapter 8.4 of \cite{torus}). It follows that $T_2(9,2,1)$ is arc-transitive, making it one of the four possible cubic arc-transitive tricirculant graphs (called $F054A$ in \cite{symtric}, following Foster's notation).

Theorem \ref{theo:t2} now follows from Lemmas \ref{lem:yk}, \ref{lem:t2cong}, \ref{lem:t2bic} and Proposition \ref{prop:torus}.
\begin{figure}[H]
\centering
\includegraphics[width=0.5\textwidth]{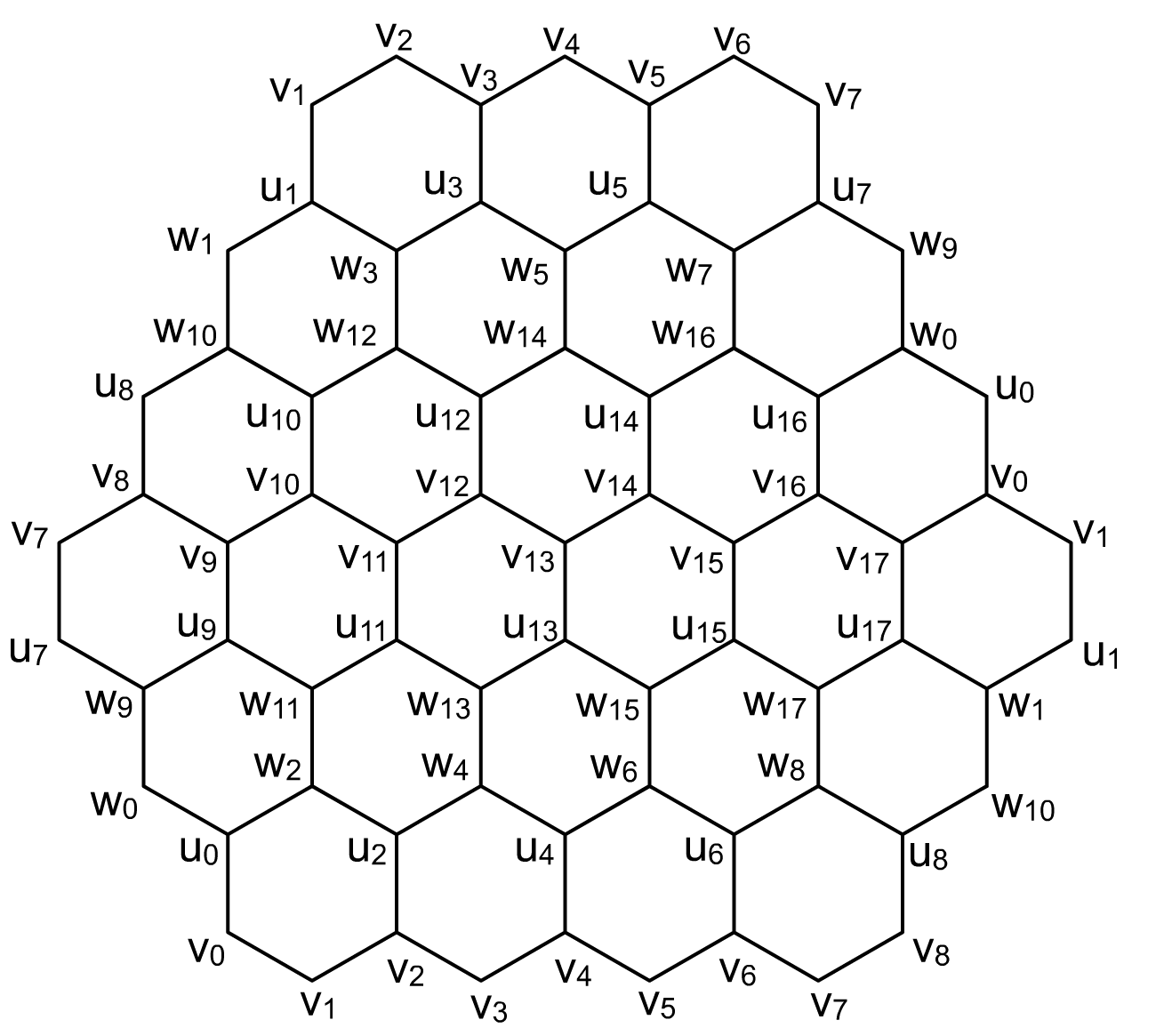}
\caption{A drawing of $T_2(9,2,1)$ as the map on the torus $\{6,3\}_{3,3}$.}
\label{fig:torus2}
\end{figure}

\section{Type 3}

\begin{figure}[H]
\centering
\includegraphics[width=0.3\textwidth]{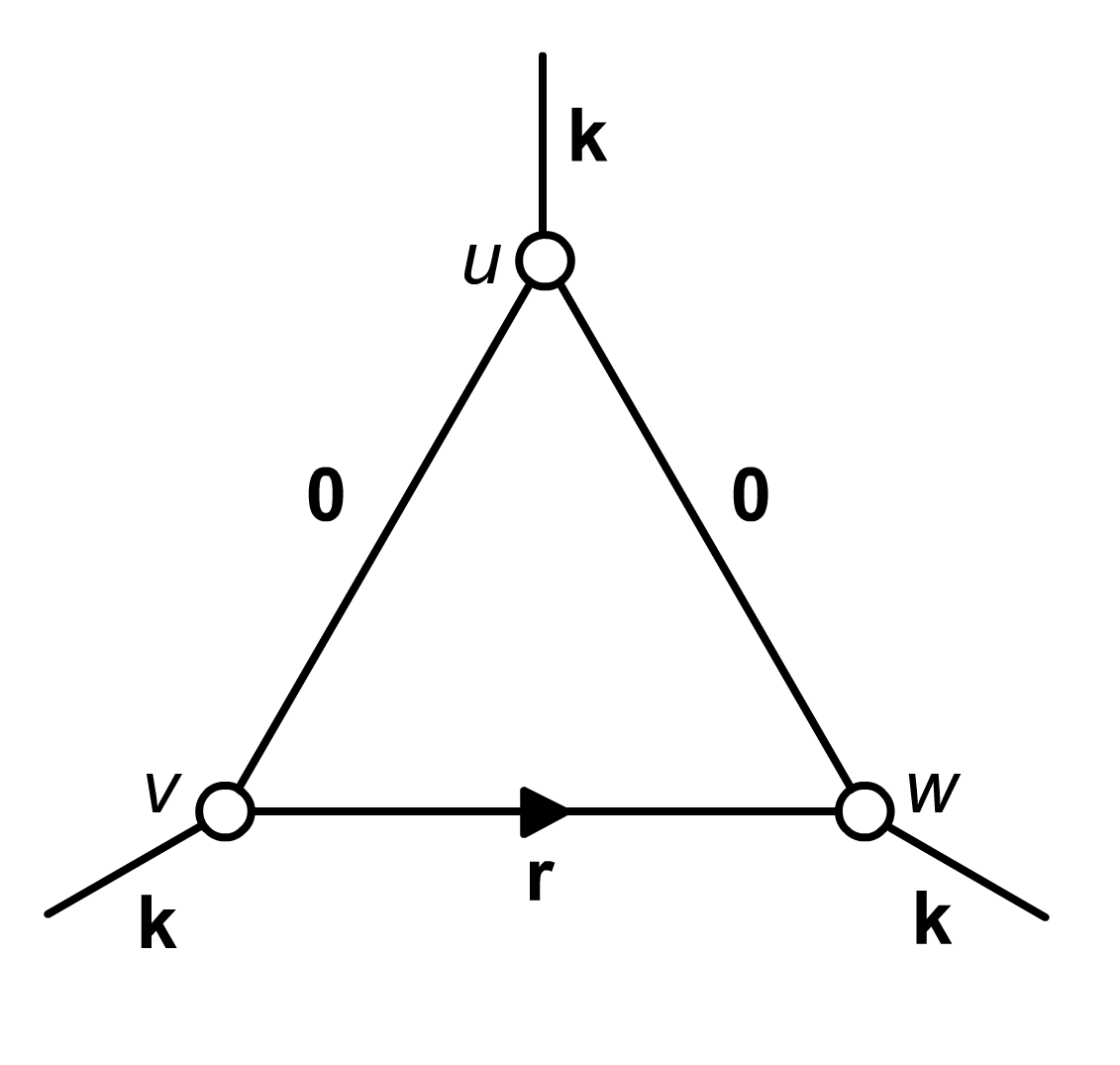}
\caption{Type 3}
\label{fig:Delta3}
\end{figure}

Let $k \geq 9$ be an integer, and let $r \in \ZZ_{2k}$. Define $T_3(k,r)$ as the derived graph of $\Delta_3$ with the normalized voltage assignment for $\ZZ_{2k}$ shown in Figure \ref{fig:Delta3}\\
Let $U=\{u_0,u_1,...,u_{2k-1}\}$, $V=\{v_0,v_1,...,v_{2k-1}\}$ and $W=\{w_0,w_1,...,w_{2k-1}\}$ be the respective fibers of vertices $u$, $v$ and $w$ in $\Delta_3$. It is clear that any cubic tricirculant of Type $3$ is isomorphic to $T_3(k,r)$ for some $k$ and $r$. 

\begin{theorem}
If $T_3(k,r)$ is connected, then it is isomorphic to either a prism or a M\"obius Ladder.  
\end{theorem}


\begin{proof}
First, recall that $T_3(k,r)$ is connected if and only if $\gcd(k,r)=1$. Then $\gcd(2k,r)~\in~\{1,2\}$, depending on whether $r$ is even or odd. 

If $r$ is odd, then the subgraph induced by $0$- and $R$-edges is a single cycle of length $6k$, $(w_0,u_0,v_0,w_r,u_r,v_r,...,w_{2k-r},u_{2k-r},v_{2k-r},w_0)$. It is straighforward to see that $K$-edges join antipodal vertices in this cycle. Hence, in this case $T_3(k,r)$ is a M\"obius Ladder.

If $r$ is even, then the graph induced by $0$- and $R$-edges is the union of two disjoint cycles of length $3k$: one consisting of all the vertices with even index, and the other consisting on those with odd index. Observe that $K$-edges connect these two cycles creating a prism.
\end{proof}

\begin{figure}[h!]
\centering
\begin{tabular}{ccc}
\includegraphics[width=0.5\textwidth]{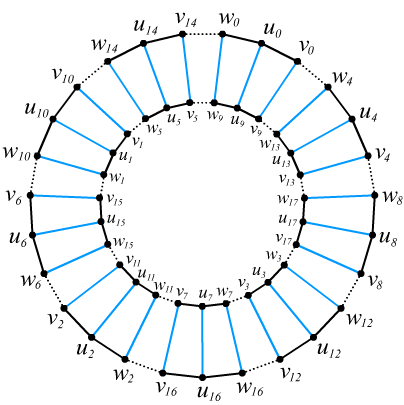} & &
\includegraphics[width=0.5\textwidth]{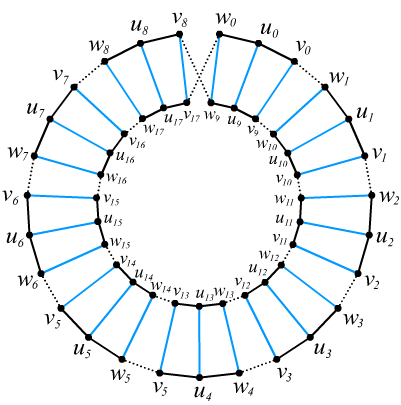} \\
$T_3(9,4)$ & & $T_3(9,1)$ 
\end{tabular}
\caption{A prism and a M\"obius Ladder on $54$ vertices.}
\label{fig:mobiusprisma}
 \end{figure}

\section{Type 4}

\begin{figure}[hhh!!!]
\centering
\includegraphics[width=0.3\textwidth]{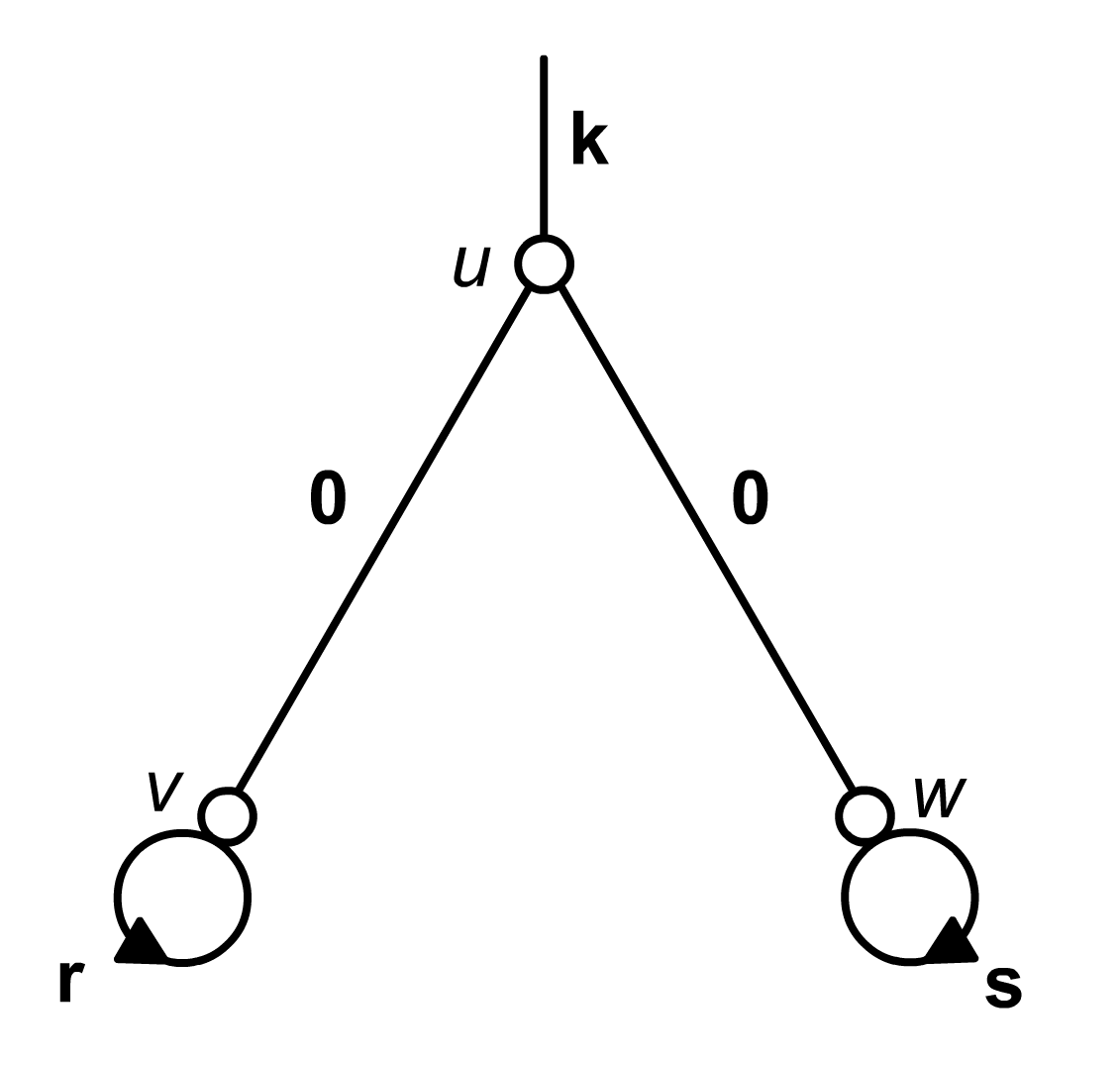}
\caption{Type 4}
\label{fig:Delta4}
 \end{figure}

Let $k$ be a positive integer, and let $r$ and $s$ be two distinct integers in $\ZZ_{2k}$. Define $T_4(k,r,s)$ as the derived graph of $\Delta_4$ with the normalized voltage assignment for $\ZZ_{2k}$ shown in Figure \ref{fig:Delta4}.
Let $U=\{u_0,u_1,...,u_{2k-1}\}$, $V=\{v_0,v_1,...,v_{2k-1}\}$ and $W=\{w_0,w_1,...,w_{2k-1}\}$ be the respective fibers of vertices $u$, $v$ and $w$ in $\Delta_4$. Then, the set of edges of $T_4(k,r,s)$ can be expressed as the union $E_K \cup E_R \cup E_S \cup E_0$ where:
\begin{eqnarray*}
E_K &=& \{u_iu_{i+k} : i \in \ZZ_{2k}\} \cr
E_0 &=& \{u_iv_i : i \in \ZZ_{2k}\} \cup \{u_iw_i : i \in \ZZ_{2k}\} \cr
E_R &=& \{w_iw_{i+r}: i \in \ZZ_{2k} \} \cr 
E_S &=& \{v_iv_{i+s}: i \in \ZZ_{2k} \} 
\end{eqnarray*}
For $X \in \{0,R,S,K\}$, edges in $E_X$ will be called edges of type $X$, or simply $X$-edges. Similarly as with tricirculants of Types $1$, $2$ and $3$, every cubic tricirculant of type $4$ with $6k$ vertices is isomorphic to $T_4(k,r,s)$ for an appropriate choice of $r$ and $s$.  

\begin{remark}
\label{rem:type4}
Note that for any $k$, $r$ and $s$ the following isomorphism holds: 
$$T_4(k,r,s) \cong T_4(k,-r,s) \cong T_4(k,r,-s) \cong T_4(k,s,r)$$
\end{remark}

\begin{theorem}
\label{theo:type4}
There are no cubic vertex-transitive tricirculants of Type $4$ of order greater or equal to $54$. 
\end{theorem} 

The rest of this section is devoted to prove Theorem \ref{theo:type4}. We will assume henceforth that $k \geq 9$ and thus the order of $T_4(k,r,s)$ is at least $54$. 

\begin{lemma}
\label{lem:rstype4}
If $T_4(k,r,s)$ is vertex-transitive, then $r \neq s$ and $r \neq -s$.
\end{lemma}

\begin{proof}
Suppose that $r=s$ and consider the graph $\Gamma := T_4(k,r,r)$, with $k \geq 9$ and $1\leq r \leq k-1$. Observe that $(v_0,v_r,u_r,w_r,w_0,u_0,v_0)$ is a $6$-cycle of $\Gamma$ that does not contain any $K$-edge but does contain edges of all other types. For $\Gamma$ to be vertex-transitive, there must be at least one $6$-cycle through a $K$-edge; otherwise $K$-edges conform a single edge-orbit and thus $U$ is a single vertex-orbit. However, such a cycle would quotient down to a closed walk of length $6$ in $\Delta_4$ having voltage $0$ and tracing the semi-edge $(uu)_k$. By observing Figure \ref{fig:Delta4} we see that
any closed walk of length $6$ visiting $(uu)_k$ has net voltage  $k \pm 3r$. This is, if $\Gamma$ is vertex-transitive, then $3r \equiv k$ \mod $2k$. Observe that under these conditions, $\Gamma$ is connected if and only if $r=1$ and thus making $k=3$, contradicting that $k \geq 9$. This shows that $r \neq s$ and view of Remark \ref{rem:type4} then $r \neq -s$.
\end{proof}

\begin{figure}[h!]
\centering
\includegraphics[width=0.6\textwidth]{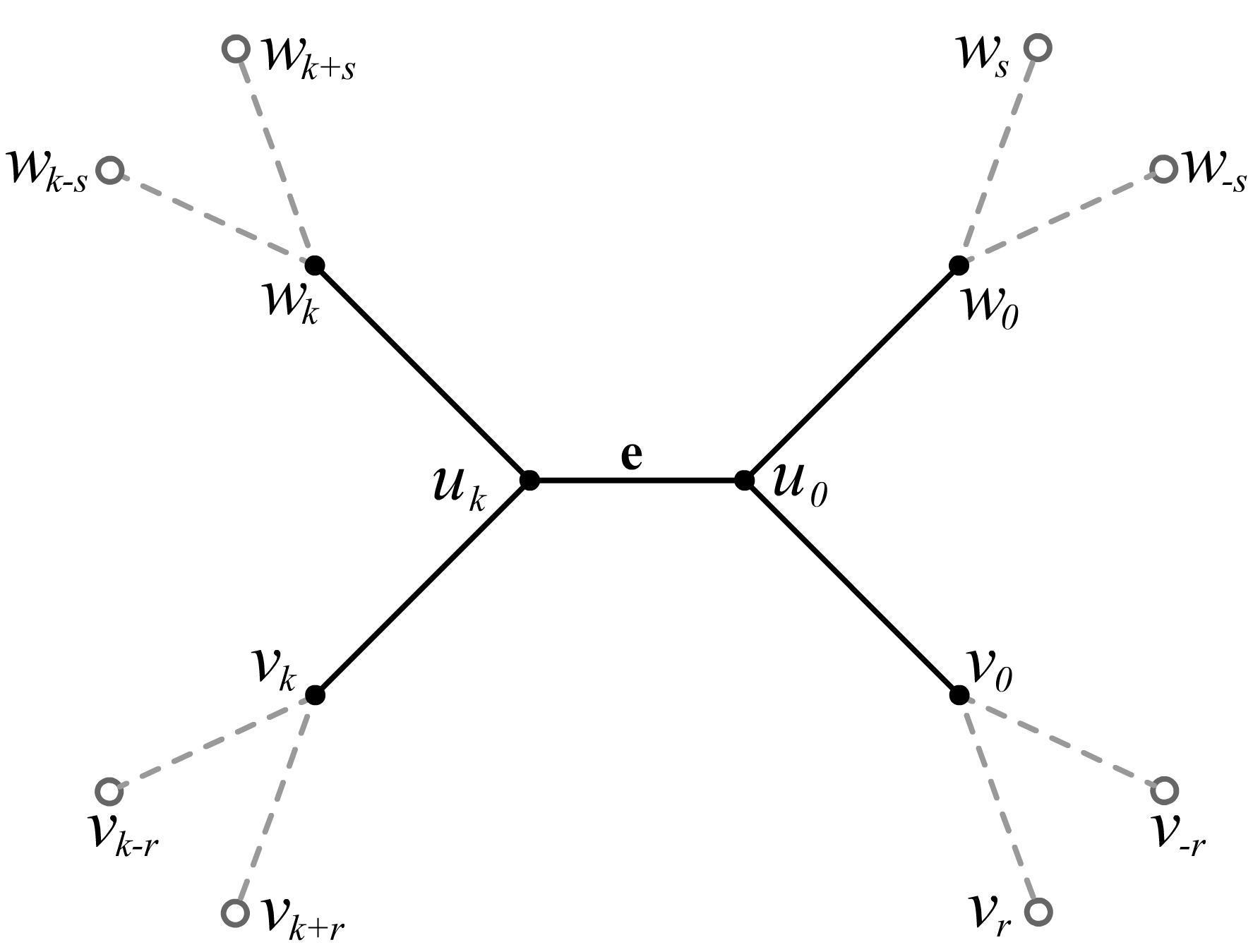}
\caption{Neighbourhood of $u_0u_k$. The subgraph $H$ shown in solid edges and vertices}
\label{fig:vecindad}
\end{figure}

\begin{lemma}
\label{lem:8cycles}
Let $e$ be a $K$-edge of $T_4(k,r,s)$ and let $a$ be an integer greater than $4$. If $C$ is an cycle of length $a$ containing $e$ then there exists another cycle of length $a$, $C'$, such that the intersection $C \cap C'$ contains a $3$-path whose middle edge is $e$.
\end{lemma}

\begin{proof}
Without loss of generality, we can assume $e=u_0u_k$. Define $H$ as the subgraph of $\Gamma$ induced by vertices at distance at most one from $e$, and let $C$ be a $k$-cycle through $e$ (see Figure \ref{fig:vecindad}). Observe that $C$ intersects $H$ in a $3$-path, $P$, whose middle edge is $e$. Now, let $\phi$ be the mapping that acts by multiplying the index of each vertex by $-1$; that is $\phi$ maps $u_i$ into $u_{-i}$, $v_i$ into $v_{-i}$, and $w_i$ to $w_{-i}$. Observe that $\phi$ is in fact an automorphism of $\Gamma$. Moreover, $\phi$ fixes each vertex and each edge of $H$. In particular, it fixes $P$. Therefore, $\phi(C)$ is a $k$-cycle through $e$ and $P$ lies in the intersection of $C$ and $\phi(C)$. To see that $\phi(C)$ is different from $C$, observe that $\phi$ interchanges the two vertices in each of the following sets: $\{v_{k+r},v_{k-r}\}$, $\{v_{r},v_{-r}\}$,$\{w_{k+s},w_{k-s}\}$ and $\{w_{s},w_{-s}\}$ (white vertices in Figure \ref{fig:vecindad}). It is clear that $C$ can visit at most one vertex in each of these four sets, and if $C$ visits one vertex in one of these sets, then $\phi(C)$ must visit the other. 
\end{proof}

\begin{lemma}
\label{lem:type41}
If $T_4(k,r,s)$ is vertex-transitive, then in $\ZZ_{2k}$ one of the equalities (A1)--(A4) below and one of the equalities (B1)--(B4) below hold:

\begin{table}[H]

\begin{tabular}{c c c c}

 (A1) & $k +2r +s = 0 \qquad \qquad$  & (B1) & $k +2s +r = 0$ \\
 (A2) & $k -2r +s = 0 \qquad \qquad$ & (B2) & $k -2s +r = 0$ \\
 (A3) & $r+3s = 0 \qquad \qquad$ & (B3) & $s+3r = 0$ \\
 (A4) & $r-3s = 0 \qquad \qquad$ & (B4) & $s-3r = 0$ \\

\end{tabular}

\label{table:type4}
\end{table}
\end{lemma}

\begin{figure}[h!]
\centering
\includegraphics[width=0.60\textwidth]{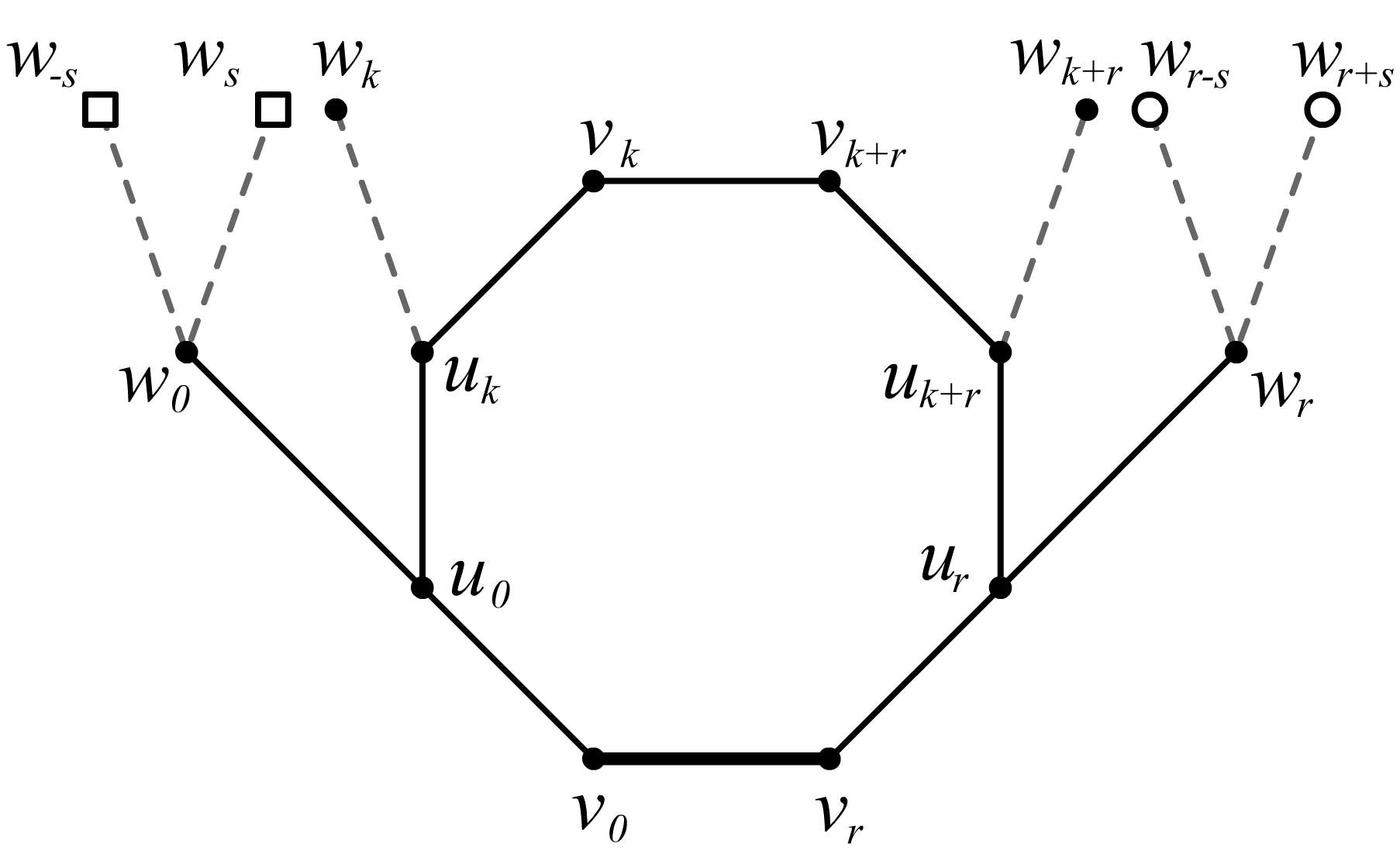}
\caption{A subgraph containing an $(v_0,v_r)$}
\label{fig:vecindadr}
\end{figure}

\begin{proof}

Suppose $T_4(k,r,s)$ is vertex-transitive. We will first show that one of the equalities (A1)--(A4) holds in $\ZZ_{2k}$. The rest of the lemma will follow from the fact that $T_4(k,r,s) \cong T_4(k,s,r)$. Since $T_4(k,r,s)$ is vertex-transitive, the edge-neighbourhood of any vertex in $U$ can be mapped by an automorphism to the edge-neighbourhood of any vertex in $V$. In particular, either there exists an automorphism mapping the $K$-edge $u_0u_k$ to the $0$-edge $u_0v_0$ or there is an automorphism mapping $u_0u_k$ to the $R$-edge $v_0v_r$. Therefore, the property described in Lemma \ref{lem:8cycles} should also hold for $u_0v_0$ or for $v_0v_r$ (or possibly both). We will see what this means in terms of $k$, $r$ and $s$.

Suppose that the property described in Lemma \ref{lem:8cycles} holds for $u_0v_0$. Observe that there is an $8$-cycle $C=(u_0,v_0,v_{-r},u_{-r},u_{k-r},v_{k-r},v_k,u_k,u_0)$ through $u_0v_0$ (see Figure \ref{fig:vecindad0}). Then there must be another $8$-cycle $C'$ whose intersection with $C$ contains the $3$-path $(v_{-r},v_0,u_0,u_k)$, but no $4$-path containing this $3$-path. This in turns implies that some vertex in $\{v_{-3r}, u_{-2r}\}$ is adjacent to a vertex in $\{w_{k-s},w_{k+s}\}$. Since no vertex in $V$ is adjacent to a vertex in $W$, we have that either $u_{-2r}w_{k-s} \in E$ or $u_{-2r}w_{k+s} \in E$. This implies $2r \equiv k+s$ or $2r \equiv k-s$, and so (A1) or (A2) holds.

\begin{figure}[h!]
\centering
\includegraphics[width=0.60\textwidth]{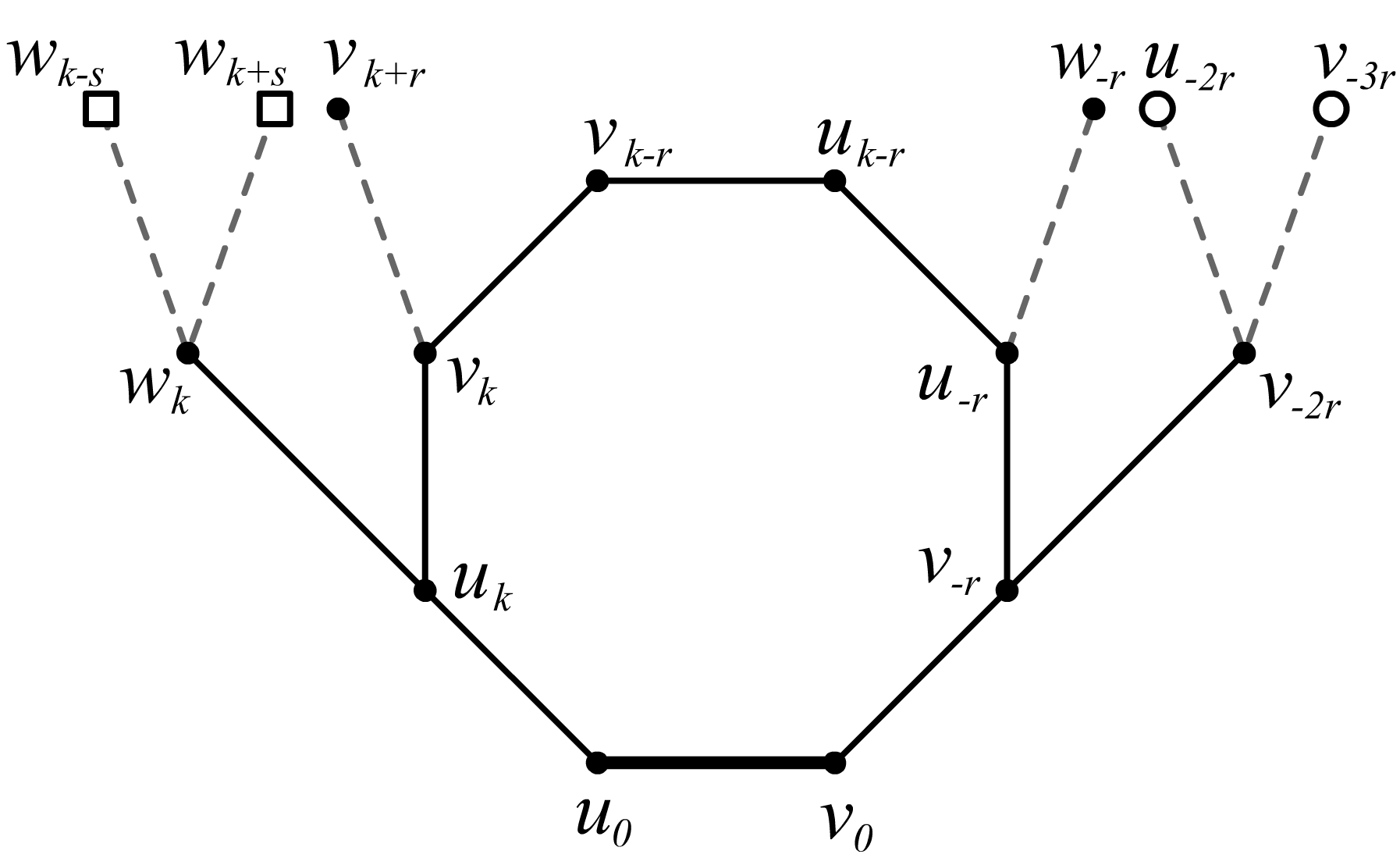}
\caption{A part of the graph $T_4(k,r,s)$ containing a $u_0v_0$}
\label{fig:vecindad0}
\end{figure}

If, on the other hand, the property described in Lemma \ref{lem:8cycles} holds for $v_0v_r$, then one vertex in $\{w_{r+s},w_{r-s}\}$ must be adjacent to a vertex in $\{w_{s},w_{-s}\}$ (see Figure \ref{fig:vecindadr}), implying that one of the following expressions must be equal to $0$ in $\ZZ_{2k}$: $r+3s$, $r+s$, $r-3s$, $r-s$. However, in view of Lemma \ref{lem:rstype4}, $r+s \neq 0$ and $r-s \neq 0$. Hence (A3) or (A4) holds. 
\end{proof}

We are now ready to prove Theorem \ref{theo:type4}. Let $k \geq 9$ be an integer and suppose $T_4(k,r,s)$ is vertex-transitive. Then by Lemma \ref{lem:type41}, one of the equalities (A1)--(A4) and one of the equalities (B1)-(B4) must hold. If, for instance both (A1) and (B1) hold, then by adding them we see that $r=-s$, which contradicts Lemma \ref{lem:rstype4}. Similarly, we get that $r \pm s$ if both equalities in any of the following pair hold: (A1) and (B2), (A2) and (B2), (A3) and (B3). If (A1) and (B4), or (A2) and (B4) hold, then $k=r$, which contradicts the simplicity of $\Gamma$. If (A1) and (B3), (A2) and (B3), or (a3) and (B4) hold, then either $k=5$ or $\gcd(k,r,s) \neq 1$. It is readily seen that we get a contradiction in each of the $16$ possible cases that arise from Lemma \ref{lem:type41}. We conclude that a connected $T_4(k,r,s)$ cannot be vertex-transitive if $k \geq 9$.
\bigskip

{\bf Acknowledgements.} The authors gratefully acknowledge support of the Slovenian Research Agency by financing the Research program P1-0294 and the second listed author Young Researcher scholarship.

\end{document}